\documentclass[12pt]{article}

\usepackage{latexsym}
\usepackage[all]{xy}
\usepackage{amsmath}
\usepackage{amssymb}
\usepackage{amsfonts}
\usepackage{color}
\usepackage{theorem}
\usepackage{mathtools}

\usepackage{hyperref}
\usepackage{xcolor}
\hypersetup{colorlinks=true,urlcolor=blue,linkcolor=blue,citecolor=blue}

\newcommand {\Vect} {\mathsf{Vect}}

\newcommand {\Map} {\mathbf{Map}}

\newcommand {\OO} {\mathcal{O}}

\newcommand {\fb} {\widehat{\partial}}
\newcommand {\A} {\mathcal{A}}
\newcommand {\M} {\mathcal{M}}
\newcommand {\T} {\mathbb{T}}
\newcommand{\C}{\mathcal{C}}

\newcommand  {\cdga}     {\mathbf{cdga}}

\newcommand  {\dSt}   {\mathbf{dSt}}

\newcommand{\s}{\infty}

\newcommand{\op}[1]{\operatorname{#1}}
\newcommand{\Xbar}{\text{{\sf Z}}}

\newcommand{\gp}{\mathbf{\sf grp}}
\newcommand{\bOO}{\mathbf{C}}
\newcommand{\Zbar}{\mathfrak{Z}}
\newcommand{\balpha}{\boldsymbol{\alpha}}

\makeatletter
\newcommand{\mquote}[1]{{\mathpalette\mqu@te{#1}}}

\newcommand{\mqu@te}[2]{%
  \sbox0{$\m@th#1\text{``}$}%
  \sbox2{$\m@th#1\text{\!\!''}$}%
  \sbox4{$\m@th#1#2$}%
  \ifdim\ht4>\dimexpr\ht0+1pt\relax
    \raisebox{\dimexpr\ht4-\height}{\box0}%
    #2%
    \raisebox{\dimexpr\ht4-\height}{\box2}%
  \else
    \box0 #2\box2
  \fi
}
\makeatother

\newcommand{\prolim}{\boldsymbol{\mathfrak{Lim}}}

\DeclareFontFamily{U}{rsf}{}
\DeclareFontShape{U}{rsf}{m}{n}{
  <5> <6> rsfs5 <7> <8> <9> rsfs7 <10->  rsfs10}{}
\DeclareMathAlphabet{\mathscr}{U}{rsf}{m}{n}
\newcommand{\mycal}[1]{\mathscr{#1}}

\newtheorem{introthm}{Theorem}
\newtheorem{thm}{Theorem}[section]
\newtheorem{prop}[thm]{Proposition}

\newtheorem{df}[thm]{Definition}
\newtheorem{cor}[thm]{Corollary}

{\theorembodyfont{\rmfamily} \newtheorem{rmk}[thm]{Remark}}
{\theorembodyfont{\rmfamily} \newtheorem{ex}[thm]{Example}}

\newtheorem{pb}[thm]{Problem}
\begin{document}

\title{\textbf{Poisson geometry  of the  moduli of \\
    local systems on smooth varieties}}  

\author{Tony Pantev and Bertrand To\"{e}n}


\date{May 2019}

\maketitle

\

\bigskip

\centerline{
\emph{\normalsize To Masaki Kashiwara on his 70th birthday}
}

\

\begin{abstract} We study the moduli of $G$-local
  systems on smooth but not necessarily proper
complex algebraic varieties. We show that, when considered as
derived algebraic stacks, they carry natural Poisson structures,
generalizing the well known case of curves. We also construct
symplectic leaves of this Poisson structure by fixing local
monodromies at infinity, and show that a new feature, called
\emph{strictness}, appears as soon as the divisor at infinity has
non-trivial crossings.
\end{abstract}

\tableofcontents

\section*{Introduction}
\addcontentsline{toc}{section}{Introduction}

For a smooth complex algebraic curve $X$ and a reductive group $G$, it
is well known that the moduli space $M_G(X)$ of $G$-local
systems\footnote{In this paper, '$G$-local systems' are
  representations of $\pi_1(X)$ into $G$} carries a canonical Poisson
structure (see \cite{fock-rosly,
  ghjw,goldman,guruprasad-rajan}). Moreover, the symplectic leaves of
this Poisson structure can be identified with moduli of $G$-local
systems having fixed conjugacy classes of monodromies at
infinity. This topological picture also has algebraic counterparts for
which local systems are replaced by flat bundles or Higgs bundles
possibly with irregular singularities, and is known to be compatible
with the comparison isomorphisms between these different incarnations
of the moduli problem (see for instance \cite{bo}).  However, as far
as the authors are aware, very little is known about the Poisson
geometry of local systems moduli on higher dimensional varieties
outside of the proper case.

The purpose of this note is to explore the moduli of $G$-local systems
on higher dimensional smooth open varieties, with a particular focus
on their Poisson geometry. For us, the results presented in this work
represent a first step towards an understanding of moduli of
local systems on higher dimensional varieties, with a long term goal
to extend Simpson's non-abelian Hodge theory to the non-proper
case. 

As a first comment, derived algebraic geometry is useful, and probably
unavoidable, for this project. Indeed, for a higher dimensional
compact oriented manifold $M$, it is known (see \cite{ptvv,ems}) that
the moduli of $G$-local systems on $M$ carries canonical symplectic
structure provided that
\begin{enumerate}
\item the moduli is considered as a \emph{derived
  algebraic stack} and not simply as a scheme or a stack.
\item the symplectic structures involve a \emph{cohomological shift}
  by $2-d$, where $d$ is the dimension of $M$.
\end{enumerate}
In this work, we consider the derived moduli stack $Loc_G(X)$ of
$G$-local systems on a complex smooth algebraic variety $X$ of
complexe dimension $d$. We establish two  principal results that can be
summarized  as follows:

\newpage

\begin{introthm}{(see Theorem~\ref{p4})}\label{ti}
\begin{enumerate}
\item The derived stack $Loc_G(X)$ carries a canonical
  $(2-2d)$-shifted Poisson structure. This recovers the usual Poisson
  structure when $d=1$.
\item The $(2-2d)$-Poisson structure on $Loc_G(X)$ has generalized
  symplectic leaves. In particular  $Loc_G(X,{\lambda_{\bullet}})$, the derived
  moduli of $G$-local systems with fixed conjugacy classes $\lambda_{i}$
  of monodromies at infinity is such a generalized symplectic leaf.
\end{enumerate}
\end{introthm}

\

Before describing the content of this work we need to add a couple of
comments concerning the previous result. The Poisson structure on
$Loc_G(X)$ will be constructed by using a very specific topological
property of smooth complex algebraic varieties, namely that their
\emph{boundary at infinity} is a compact manifold (of real dimension
$2d-1$ if $X$ is of dimension $d$). As a consequence, there is a
natural map $Loc_G(X) \rightarrow Loc_G(\partial X)$ sending a
$G$-local system on $X$ to its restriction to the boundary. By the
work of Calaque \cite{ca} it is known that such restriction maps come
equiped with a canonical Lagrangian structure, and by the work of
Melani and Safronov \cite{mesa1,mesa2} it is known that Lagrangian
structures induce Poisson structures. This roughly explains why
statement $(1)$ is true.  Statement $(2)$ is subtler, mainly because
one has to make precise what 'fixing the monodromies at infinity'
means. This is particularly important in the derived setting where the
fixing of the local monodromies involves higher homotopy coherence
conditions. Moreover, we only prove $(2)$ under the restrictive
condition that the divisor at infinity of $X$ can be chosen to have at
most double intersections.  We will see that even in this simple case
a new feature appears, and that we have to impose an additional
condition on the local monodromies at infinity that we call
\emph{strictness} (see Definition~\ref{d3}). This condition is
invisible on the non-derived moduli space, but is required in order to
construct symplectic leaves in the full derived moduli stack.

The paper is organized as follows.  In section~\ref{sec:1} we start
with a short reminder on the derived moduli of $G$-local systems on a
space, and the various ways in which one could describe this derived
moduli in concrete algebraic terms. In sections~\ref{sec:2} and
\ref{sec:3} we briefly recall shifted symplectic and Poisson
structures, and introduce the notion of \emph{generalized symplectic
  leaves} in this context. In section~\ref{sec:4} we focus on the case
of complex smooth algebraic varieties. We examine their structure at
infinity and deduce the existence of the shifted Poisson structure on
the derived moduli of local systems. We first analyze the special case
of a smooth divisor at infinity and show that the construction works
in essentially the same manner as in the case of curves. Finally we
study the case of a divisor with two smooth intersecting components
and show how the strictness condition appears naturally when one tries
to construct symplectic leaves. We also provide families of examples
of strict pairs. In section~\ref{sec:5} we collect some ideas
indicating how the statements of this paper can be generalized to the
de Rham setting in which local systems are replaced by bundles with
flat connections.

\

\noindent
{\bfseries Acknowledgements:} We would like to thank Sasha Efimov,
Dmitry Kaledin, Takuro Mochizuki, and Gabriele Vessozi for several
illuminating discussions on the subject of this work.

During the preparation of this work Bertrand To\"en was partially
supported by ERC-2016-ADG-741501.  Tony Pantev was partially supported by
NSF research grant DMS-1601438, by Simons Collaboration grant \#
347070, and by the Laboratory of Mirror Symmetry NRU HSE, RF
Government grant, ag. No 14.641.31.0001.
    
\subsection*{Notation and conventions} 
\addcontentsline{toc}{subsection}{Notation and conventions}

\begin{description}
\item[$k$] - a field of characteristic zero.
\item[$\cdga_{k}^{\leq 0}$] - the $\s$-category of non-positively
  graded commutative dg-algebras over $k$.
\item[$\T$] - the $\s$-category of spaces, or equivalently the
  $\s$-category of simplicial sets.
\item[$\op{Pro}(\T)$] - the $\s$-category of pro-simplicial
  sets.
\item[$G$] - a reductive group over $k$.
\item[{$[G/G]$}] - the stack quotient of $G$ by the conjugation action
  of $G$ on itself.
\item[$\lambda$] - an element in $G$.
\item[$\bOO_{\lambda}$] - the conjugacy class of $\lambda$ in $G$.
\item[{${}_{\alpha} \widetilde{[G/G]} \to D$}] - locally constant stack with
  fiber $[G/G]$ obtained by twisting the constant stack with a circle bundle
  classified by $\alpha : D \to BS^{1}$.  
\item[$\Gamma$] - a finitely presentable discrete group. 
\item[$R_{G}(\Gamma)$] - the affine $k$-scheme parametrizing group
  homomorphisms $\Gamma \to G$.
\item[$M_{G}(X)$] - the $G$-character scheme of $X$.
\item[$\mathcal{M}_{G}(X)$] - the stack of $G$-local systems on $X$.
\item[$Loc_{G}(X)$] - the derived stack of $G$-local systems on $X$.
\item[$Loc_{G}(X,\{\lambda_{i}\})$] - the derived stack of $G$-local
  systems with fixed conjugacy classes of monodromies at infinity.
\item[{$\Gamma(F,\op{Sym}_{\OO}(\mathbb{T}_{F}[-n-1]))[n+1]$}] - the
  complex of $n$-shifted polyvectors on a derived  Artin stack $F$. 
\end{description}

\section{The moduli of local systems as a derived stack} \label{sec:1}

In this section we review the basic constructions of character
schemes, the stack of local systems, as well as the derived stack of
local systems, associated to a connected finite CW complex $X$.  We
explain how to understand the derived structure on the moduli stack of
local systems by means of a \emph{free resolutions of the space
  $X$}. We also discuss the basics of differential calculus on this
derived stack by presenting an explicit model for computing algebraic
de Rham cohomology. Most of the material here is well known or at
least part of the folklore.

\subsection{The character scheme
  and the stack of local systems}

Let $k$ be a field of characteristic zero, $X$ a connected finite CW
complex, and $G$ a reductive group over $k$.  We consider $G$-local
systems on $X$ which are by definition locally constant principal
$G$-bundles on $X$. If we fix a base point $x \in X$ we can
equivalently view $G$-local systems as $G$-valued representations of
the discrete group $\Gamma:=\pi_{1}(X,x)$.

The \emph{\bfseries moduli of $G$-local systems}
can then be defined by 
$$
M_{G}(X):=R_{G}(\Gamma)/G:=\op{Hom}_{\gp}(\Gamma,G)/G.
$$
This formula can have several interpretations, depending on how we
view its terms. In the most straightforward interpretation (see
e.g. \cite{lubotzky-magid}) $R_{G}(\Gamma)=\op{Hom}_{\gp}(\Gamma,G)$
is an affine scheme over $k$, classifying group homomorphisms $\Gamma
\rightarrow G$. It can be constructed explicitly as a closed subscheme
in $G^{p}$ where $p$ is the number of a chosen set of generators for
$\Gamma$ and the ideal cutting out $R_{G}(\Gamma)$ is given by the
relations among these generators. Alternatively we can define
$R_{G}(\Gamma)$ as the affine $k$-scheme which represents the functor
sending a commutative $k$-algebra $A$ to the set
$\op{Hom}_{\gp}(\Gamma,G(A))$ of group homomorphisms from $\Gamma$ to
the group of $A$-points $G(A)$.  The group $G$ acts on $R_{G}(\Gamma)$
by conjugation.  The quotient of $R_{G}(\Gamma)$ by $G$ can itself be
interpreted as an affine GIT quotient. Thus $M_{G}(X)$ is an affine
scheme over $k$ whose ring of functions is the ring of $G$-invariant
functions on $R_G(\Gamma)$. The set of $k$-points of $M_{G}(X)$ is in
one-to-one correspondence with the set of isomorphism classes of
semi-simple locally constant principal $G$-bundles on $X$. The scheme
$M_{G}(X)$ is also often called the $G$-character scheme of $X$.

A less naive viewpoint is to consider the quotient stack of
$R_{G}(\Gamma)$ by the action of the group $G$. This stack, denoted by
$\M_{G}(X)=[R_G(\Gamma)/G]$, is called the \emph{\bfseries stack of
  $G$-local systems} on $X$. The $k$-points of $\M_G(X)$ form a
groupoid equivalent to the groupoid of all $G$-local systems on
$X$. The stack $\M_{G}(X)$ is an algebraic stack in the sense of Artin
and comes equipped with a structure morphism $\M_{G}(X)
\longrightarrow M_{G}(X)$ which is universal for morphism to
schemes. In other words the character variety $M_{G}(X)$ is a coarse
moduli space for the stack $\M_{G}(X)$.

\subsection{Simplicial resolutions and the derived stack of local
  systems}

In this work we will need a slightly more refined version of the stack
of local systems called the \emph{\bfseries derived stack of local
  systems}. The derived stack of local systems arises naturally both
as a way of encoding the algebraic complexity of the relations
defining $\Gamma$ and as a device for repairing singularities in
$\M_{G}(X)$.

The scheme $R_{G}(\Gamma)$ and hence the stack $\M_{G}(X)$ can in
general be very singular. However, when the group $\Gamma$ happens to
be free of rank $p$, $R_{G}(\Gamma)$ is isomorphic to $G^{p}$ and is
thus smooth over $k$. When $\Gamma$ is not free we can consider
\cite[Chapter VI]{may}, \cite[Chapter 5]{goerss-jardine} a simplicial
free resolution $B\Gamma_{\bullet} \simeq X$ of the space $X$. More
precisely, given a base point $x \in X$, consider a simplicial group
model for the loop group $\Omega_{x}(X)$. This simplicial group can be
resolved by free groups, i.e. replaced by a weakly equivalent
simplicial group $\boldsymbol{\Gamma}_{\bullet}$ where each
$\boldsymbol{\Gamma}_{n}$ is free on a finite number of generators.
Note that the geometric realization of the simplicial space
$B\boldsymbol{\Gamma}_{\bullet}$ is homotopy equivalent to $X$, and we
can thus view $\boldsymbol{\Gamma}_{\bullet}$ as a free resolution of
the pointed space $(X,x)$. Note also that this resolution depends on $X$
and not just on the group $\Gamma$ (except when $X$ is itself a
$K(\Gamma,1)$ in which case $\boldsymbol{\Gamma}_{\bullet}$ is a free
resolution of the group $\Gamma$).  Applying $R_{G}(-)$ to
$\boldsymbol{\Gamma}_{\bullet}$ yields a cosimplicial affine scheme
$R_{G}(\boldsymbol{\Gamma}_{\bullet})$, or equivalently a simplicial
commutative $k$-algebra
$\OO(R_{G}(\boldsymbol{\Gamma}_{\bullet}))$. The passage to normalized
chains defines a commutative dg-algebra, whose quasi-isomorphism type
does not depend on the choice of the resolution
$\boldsymbol{\Gamma}_{\bullet}$ of $X$. In other words we get a
commutative dg-algebra $\mycal{A}_{G}(X)$ which, up to quasi-isomorphism,
only depends on the homotopy type of $X$.

By construction $H^{0}(\mycal{A}_{G}(X))$ is naturally isomorphic to
$\OO(R_{G}(\Gamma))$ and $H^{i}(\mycal{A}_{G}(X))$ vanish for $i>0$. The
other cohomologies $H^{i}(\mycal{A}_{G}(X))$ for $i<0$ can be non-zero.  When
$X=K(\Gamma,1)$, the cohomology $H^{\bullet}(\mycal{A}_{G}(X))$ is the so
called \emph{\bfseries representation homology} of the group ring
$k[\Gamma]$ in the sense of \cite{bkr} and codifies many interesting
invariants of the group $\Gamma$. For an arbitrary CW complex $X$ the
$k$-vector spaces $H^{i}(\mycal{A}_{G}(X))$ are invariants of the space $X$
and may be non-trivial even when $X$ is simply connected (see
Example~\ref{celldec}).

As explained in \cite{ems} the non-positively graded cdga $\mycal{A}_G(X)$
has a spectrum $Spec\,\mycal{A}_{G}(X)$ which is a \emph{\bfseries
  derived affine scheme}, that is  -  an affine $k$-scheme equipped with a
sheaf of cdga. The conjugation action of $G$ on the various
$R_{G}(\boldsymbol{\Gamma}_{n})$ gives rise to an action on the
commutative dg-algebra $\mycal{A}_{G}(X)$ and hence $G$ acts on its
spectrum. The quotient stack
\[
  Loc_{G}(X) := [Spec\,\mycal{A}_{G}(X)/G]
\]
is the
\emph{\bfseries derived stack of $G$-local systems on $X$} of
\cite{hagdag}.  We refer the reader to \cite{hagII} for the formalism
of derived schemes and derived stacks, in particular we will not
explain in this work how to formally construct the $\s$-category of
derived stacks and how to define the above quotient.

Note that, as explained in \cite{hagdag}, 
$Loc_{G}(X)$ can also be considered as an ($\infty-$)functor
$$
Loc_{G}(X) : \cdga_k^{\leq 0} \longrightarrow \T
$$
on the $\s$-category $\cdga_{k}^{\leq 0}$ of non-positively graded
commutative $k$-linear dg-algebras.  This functor sends a dg-algebra
$A$ to the simplicial set $Map(S(X),BG(A))$ of maps from the singular
simplicies in $X$ to the simplicial set of $A$-points of the stack
$BG$ (see \cite{hagdag} for details). In the special case when $G=
GL_{n}$ the simplicial set $Loc_{G}(X)(A)$ also admits an alternative
sheaf theoretic description.  Consider the category whose objects are
sheaves of $A$-dg-modules on $X$ that are locally quasi-isomorphic to
the constant sheaf $A^{\oplus n}$, and whose morhphisms are
quasi-isomorphisms between such sheaves. The nerve of this category is
naturally equivalent to $Loc_{G}(X)(A)$ (see \cite{hagdag}).

\begin{ex}\label{celldec} There is  another useful description of
the derived stack $Loc_G(X)$ which instead of a free resolution
$\boldsymbol{\Gamma}_{\bullet}$ of the space $X$ uses a cell
decomposition of $X$ as follows.

Let us assume that we have fixed a cell decomposition of $X$:
\[
\xymatrix@1@M+0.5pc{\varnothing = X_0 \ar@{^{(}->}[r] & \cdots
  \ar@{^{(}->}[r] & X_k \ar@{^{(}->}[r] & X_{k+1} \ar@{^{(}->}[r] &
  \cdots \ar@{^{(}->}[r] & X_n=X,}
\]
where each inclusion $X_k \hookrightarrow X_{k+1}$ is obtained by a
push-out
$$
\xymatrix{
X_k \ar@{^{(}->}[r] & X_{k+1} \\
S^{n_k} \ar@{^{(}->}[r] \ar[u] & B^{n_k+1} \ar[u]
}
$$
adding a $(n_k+1)$-dimensional cell. The derived stack $Loc_G(X)$
itself decomposes as a tower of maps
\[
\xymatrix@1@M+0.5pc{
  Loc_G(X) \ar[r] & \cdots \ar[r] & Loc_G(X_{k+1}) \ar[r] &
  Loc_G(X_k) \ar[r] & \cdots \ar[r] & Loc_G(\varnothing)=*,}
\]
where each map is a part of a pull-back square
$$
\xymatrix{
Loc_G(X_{k+1}) \ar[d] \ar[r] & Loc_G(X_k) \ar[d] \\
Loc_G(B^{n_k+1})=BG \ar[r] & Loc_G(S^{n_k}).
}
$$
Moreover, for an $m$-dimensional sphere $S^m$, the derived stack
$Loc_G(S^m)$ can be computed explicitly, for instance by induction on
$m$ using the cell decomposition $S^m=B^{m}\sqcup_{S^{m-1}}B^m$. We
have
$$Loc_G(S^0)\simeq BG \times BG \qquad Loc_G(S^1) \simeq [G/G]$$
and for any $m>1$ 
$$Loc_G(S^m)\simeq [Spec\, \mycal{A}_G(S^m)/G]$$
with 
$\mycal{A}_G(S^m)\simeq Sym_k(\frak{g}^\vee[m-1])$, and where 
$\frak{g}^\vee$ is the $k$-linear dual of the Lie algebra of $G$.

\end{ex}

\subsection{Cotangent complexes and differential forms}

The derived stack $Loc_{G}(X)$, being a quotient of a derived affine
scheme by an algebraic group, is a derived Artin stack. In particular
it has a cotangent complex $\mathbb{L}_{Loc_{G}(X)}$ which is a
quasi-coherent complex on $Loc_{G}(X)$. This can be described
explicitly in terms of the $G$-equivariant dg-algebra
$\mycal{A}_{G}(X)$ as follows.

First note that  the $\s$-category of quasi-coherent complexes on
$[Spec\, \mycal{A}_{G}(X)/G]$ is naturally equivalent to the $\s$-category of
$G$-equivariant $\A_{G}(X)$-dg-modules. The derivative of the $G$-action
on $\mycal{A}_{G}(X)$ induces a morphism of $\mycal{A}_{G}(X)$-dg-modules
$$
a : \mathbb{L} \longrightarrow \mathfrak{g}^{\vee}\otimes_{k} \mycal{A}_{G}(X),
$$
where $\mathbb{L}$ is the cotangent complex of the commutative
dg-algebra $\mycal{A}_{G}(X)$, and $\mathfrak{g}^{\vee}$ is the dual of the
Lie algebra of $G$. The homotopy fiber of $a$ is a well defined
$\mycal{A}_{G}(X)$-module which, considered as a quasi-coherent module on
$Spec\, \mycal{A}_G(X)$, is naturally equivalent to the pull-back of the
cotangent complex of $Loc_G(X)$ by the atlas map $Spec\, \mycal{A}_G(X)
\rightarrow Loc_G(X)$. This homotopy fiber carries a natural
\linebreak
$G$-equivariant structure, making it into a quasi-coherent module on
$Loc_G(X)$, which is the cotangent complex of $Loc_G(X)$.

Global sections are easier to understand in this setting, and are
simply obtained by taking $G$-invariants:
$$
a^G : \mathbb{L}^G \longrightarrow
(\mathfrak{g}^{\vee}\otimes_{k} \mycal{A}_{G}(X))^G.
$$
The homotopy fiber of the map $a^G$ computes
$\Gamma(Loc_G(X),\mathbb{L}_{Loc_G(X)})$, the complex of global
sections of the cotangent complex.

More generally, as explained in \cite{ptvv,ems}, we can talk about
differential forms and the whole de Rham complex (endowed with its
natural Hodge filtration) on the derived Artin stack
$Loc_{G}(X)$. This is a complex $\A^{\bullet}(Loc_{G}(X))$, filtered
by subcomplexes
$F^p\A^{\bullet}(Loc_{G}(X)) \subset \A^{\bullet}(Loc_{G}(X))$. The
complex $\A^{\bullet}(Loc_{G}(X))$ computes de Rham cohomology of
$Loc_{G}(X)$, while  the complex
\[
\A^{p,cl}(Loc_{G}(X)):=F^p A^{\bullet}(Loc_{G}(X))[p]
\]
is called the
\emph{\bfseries complex of closed $p$-forms on $Loc_{G}(X)$}.  In our
 setting, these complexes can be described explicitly as follows.  With
the same notations as above, we form the graded $k$-module
$$
C:=(\op{Sym}_{\mycal{A}_{G}(X)}(\mathbb{L}[-1])\otimes_k
\op{Sym}_k(\mathfrak{g}^{\vee}[-2]))^G.
$$
This graded module comes equipped with a differential which is the sum
of three different terms: the internal cohomological differential of
$\mycal{A}_{G}(X)$, the differential induced by the coaction map \linebreak
$a: \mathbb{L} \rightarrow \mathfrak{g}^{\vee}\otimes_k \mycal{A}_{G}(X)$, and
the de Rham differential on
$\op{Sym}_{\mycal{A}_{G}(X)}(\mathbb{L}[-1])$. This makes $C$ into a complex
of $k$-modules. Moreover, $C$ comes equipped with a natural Hodge
filtration, i.e. the stupid filtration for the natural grading on
$ \op{Sym}_{\mycal{A}_{G}(X)}(\mathbb{L}[-1])\otimes_k
\op{Sym}_k(\mathfrak{g}^{\vee}[-2]))$.  The complex $C$ with this
filtration is a model for the filtered complex
$\A^{\bullet}(Loc_{G}(X))$.

\section{Symplectic and Lagrangian structures} \label{sec:2}

Recall from \cite{ptvv,cptvv} the notions of shifted symplectic and
Poisson structures on derived Artin stacks. As we noted above, for a
derived Artin stack $F$ we have a complex of closed $2$-forms
$\A^{2,cl}(F)$, defined as the second layer in the Hodge filtration on
its de Rham complex (shifted by $2$).  An $n$-cocyle in the complex
$\A^{2,cl}(F)$ is called a \emph{\bfseries closed $2$-form of degree
  $n$} (see \cite{ptvv}).  Such a form is furthermore \emph{\bfseries
  non-degenerate} if the contraction with the induced element in
$H^n(F,\wedge \mathbb{L}_{F})=H^n(\A^{2,cl}(F))$ gives a
quasi-isomorphism of quasi-coherent complexes $\omega^{\flat} :
\mathbb{T}_{F} \; \widetilde{\to} \; \mathbb{L}_{F}[n]$. A
non-degenerate closed $2$-form of degree $n$ on $F$ is called an
\emph{\bfseries $n$-shifted symplectic structure}.  This notion of
symplectic structure can be extended to the relative setting and gives
rise to the notion of a Lagrangian structure. For a morphism $f : F
\rightarrow F'$ between derived Artin stacks, an \emph{\bfseries
  $(n-1)$-shifted isotropic structure on $f$} consists by definition
of a pair $ (\omega,h)$, where $\omega$ is an $n$-shifted symplectic
structure on $F'$, and $h$ is a homotopy between $f^*(\omega)$ and $0$
inside the complex $\A^{2,cl}(F)$, i.e. $h$ is a degree $(n-1)$
cochain in $\A^{2,cl}(X)$ with coboundary $f^*(\omega)$. Such an
isotropic structure is a \emph{\bfseries $(n-1)$-shifted Lagrangian
  structure} if moreover the induced canonical morphism $h^{\flat} :
\mathbb{T}_{f} \; \widetilde{\to} \; \mathbb{L}_{F}[n-1]$ from the
relative tangent complex $\mathbb{T}_f$ of $f$ to the shifted
cotangent complex of $F$ is a quasi-isomorphism.

As shown in \cite{ptvv}, when $X$ is a compact oriented manifold of
dimension $d$, the derived stack $Loc_{G}(X)$ has a natural
$(2-d)$-shifted symplectic structure. This structure is canonical up
to a choice of a non-degenerate element in
$(\op{Sym}^{2}\mathfrak{g}^{\vee})^{G}$ which always exists since $G$
is assumed to be reductive. This statement can be extended to a
compact oriented manifold $X$ with non-empty boundary $\partial X$. By
\cite{ca}, the induced restriction map
$$
Loc_{G}(X) \longrightarrow Loc_{G}(\partial X)
$$ carries a canonical $(2-d)$-shifted Lagrangian structure for which
the $3-d = 2-(d-1)$-shifted \linebreak symplectic structure on the
target is the one discussed above.  When $\partial X=\varnothing$ we
have that \linebreak $Loc_{G}(\partial X) = Loc_{G}(\varnothing) = *$
is a point and the Lagrangian structure on $Loc_{G}(X) \rightarrow *$
recovers the $(2-d)$-shifted symplectic structure on $Loc_{G}(X)$.

In fact, in order to get a Lagrangian structure on a map between
moduli of local systems it is not necessary for the map to be induced
from restricting to an actual boundary. Indeed, for a continuous map
between finite CW complexes $f : Y \rightarrow X$, there is a notion
of an \emph{\bfseries orientation of dimension $d$} on $f$. By
definition such an orientation is given by of a morphism of complexes
$\mathsf{or} : C^{\bullet}(Y,X) \longrightarrow k[1-d]$, where
$C^{\bullet}(Y,X)$ is the cofiber of the pull-back map
$f^*C^{\bullet}(X) \rightarrow C^{\bullet}(Y)$ on singular cochains
with coefficients in $k$. The morphism $\mathsf{or}$ is also assumed
to satisfy a non-degeneracy condition that ensures Poincar\'e duality
between $H^*(X)$ and $H^*(X,Y)$. Concretely, the composition of
the cup-product on
$C^{\bullet}(X)$ with the orientation map produces a well defined
pairing
$$
C^{\bullet}(X) \otimes C^{\bullet}(X,Y) \longrightarrow k[1-d]
$$
and we require that this pairing is non-degenerate on cohomology
and induces a quasi-isomorphism $C^{\bullet}(Y,X) \simeq C^{\bullet}(X)^*[1-d]$.

By \cite{ca}, when $f: Y \rightarrow X$ is endowed with an orientation
of dimension $d$, the pullback map on the derived stacks of local
systems $f^{*} : Loc_{G}(Y) \longrightarrow Loc_{G}(X)$ carries a
canonical $(2-d)$-shifted Lagrangian structure (again up to a choice of
a non-degenerate element in $\op{Sym}^{2}(\mathfrak{g}^{\vee})^{G}$).

\begin{ex} \label{ex:surface} In the special case where $X$ is a
  Riemann surface with boundary $\partial X$, we expect this to match
  the well known symplectic structures on moduli of $G$-local systems
  on $X$ with prescribed monodromies at infinity which are usually
  constructed by quasi-Hamiltonian reduction (see
  \cite{almame}). Indeed, here $Y = \partial X$ is a disjoint union of
  oriented circles, and we thus have $Loc_{G}(Y)\simeq \prod \ [G/G]$
  where $[G/G]$ denotes the stack quotient of $G$ by its conjugation
  self action. The stack $Loc_{G}(S^1)=[G/G]$ carries a canonical
  symplectic structure of degree $1$. Moreover, for any element
  $\lambda \in G$ with conjugacy class $\bOO_{\lambda}$, and
  centralizer $G_{\lambda}$ the inclusion $\bOO_{\lambda} \subset G$
  produces a canonical Lagrangian structure $BG_{\lambda}\simeq
  [\bOO_{\lambda}/G] \subset [G/G]$. Therefore, by choosing a family
  of elements $\lambda_{i} \in G$, we have two $0$-shifted
  Lagrangian morphisms
$$
\xymatrix{ \prod BG_{\lambda_{i}} \ar[rd] & &
  Loc_{G}(X). \ar[ld] \\
& \prod \ [G/G] & }
$$
By \cite{ptvv} the fiber product of these two maps therefore comes
equipped with a $0$-shifted symplectic structure. This fiber product,
denoted by $Loc_{G}(X,\{\lambda_{i}\})$ is the derived stack of
$G$-local systems on $X$ whose local monodromies at infinity are
required to belong to the conjugacy classes $\bOO_{\lambda_{i}}$. This
should recover the symplectic structures of
\cite{almame,fock-rosly,ghjw,goldman,guruprasad-rajan}. However, the
precise analysis of this comparison has not been carried out in the
literature.
\end{ex}

\section{Poisson structure and generalized
  symplectic leaves} \label{sec:3}

There is a notion of a shifted Poisson structure, generalizing the
notion of a shifted symplectic structure.  The definitions can be
found in \cite{cptvv,pridham-poisson,pv} but are long and technical
and will not be discussed here. We only recall that for any derived
Artin stack $F$ we can form the complex of \emph{\bfseries
  $n$-shifted polyvectors}
$\Gamma(F,\op{Sym}_{\OO}(\mathbb{T}_F[-n-1]))[n+1]$ which carries a
canonical Lie bracket making it into a graded dg-Lie algebra (see
\cite{cptvv} for details).  By definition, an \emph{\bfseries
  $n$-shifted Poisson
structure} on $F$ consists of a morphism in the $\s$-category of graded
dg-Lie algebras
$$
p : k[-1](2) \longrightarrow
\Gamma(F,\op{Sym}_{\OO}(\mathbb{T}_F[-n-1]))[n+1],
$$
where $k[-1](2)$ is the graded dg-Lie algebra which is just $k$ placed
in homological degree $1$ and grading degree $2$, equipped with the
zero Lie bracket.

One of the main comparison results of \cite{cptvv,pridham-poisson}
states that the space of all $n$-shifted symplectic structures on a
derived stack $F$ is equivalent to the space of all non-degenerate
$n$-shifted Poisson structures. This result was recently generalized
\cite{mesa1,mesa2} to Lagrangian structures: the space of all
Lagrangian structures on a morphism $F \rightarrow F'$ is equivalent
to the space of all non-degenerate coisotropic structures. In
particular, an $n$-shifted Lagrangian morphism of derived Artin stacks
$F \rightarrow F'$ always induces an $n$-shifted Poisson structure on
$F$. Moreover, it is expected that all $n$-shifted Poisson structures
arise this way as soon as one allows $F'$ to be a formal derived stack
rather than a derived Artin stack. Here the formal derived stack $F'$ is
the quotient of $F$ by a derived Lie algebroid defined by the
Poisson structure. Although the quotient $F'$ can be given a precise
meaning (see for instance \cite{nu}), the fact this can be enhanced to
an equivalence between shifted Poisson structures on $F$ and
Lagrangian maps out of $F$ has been announced by Costello-Rozenblyum
but has not been written yet. At any rate, this suggests that we can
view an $n$-shifted Poisson structure on a given derived stack $F$ as
an \emph{\bfseries equivalence class of $n$-shifted Lagrangian
  morphisms $F \rightarrow F'$} with $F'$ possibly formal derived
Artin stack. Here two such morphisms $F \rightarrow F'$ and
$F\rightarrow F''$ are declared to be equivalent if there is a third
one $F \rightarrow G$ and a commutative diagram
$$
\xymatrix{
  &  F' \\
  F  \ar[r] \ar[ru] \ar[rd] & G \ar[u]_-{a} \ar[d]^-{b} \\
  & F''}
$$
with $a$ and $b$ formally \'{e}tale and compatible with the Lagrangian
structures.  In fact any morphism $f : F \rightarrow F'$ can be
factored as $\xymatrix@1{ F \ar[r] & \widehat{F} \ar[r] & F'}$ where
$\widehat{F}$ is the formal completion of the morphism $f$ and
$\widehat{F} \rightarrow F$ is \'etale. The derived stack
$\widehat{F}$ is only a formal stack in general, and can be obtained
as the quotient of $F$ by the action of the Lie algebroid induced from
the morphism $f$, i.e. the Lie algebroid corresponding to the relative
tangent complex $\mathbb{T}_{f}$ of $f$. Intuitively this quotient
contracts infinitesimally all fibers of $f$. Therefore, the derived
stack $G$ in the above diagram can always be taken to be $\widehat{F}$
(for one of the two morphisms).  Conversely, given an $n$-shifted
Poisson structure on $F$, one can define from it a symplectic Lie
algebroid in the sense of \cite{pysa}, whose quotient is expected to
recover a Lagrangian map $F \longrightarrow F'$ that induces back the
original Poisson structure on $F$.

Even though this point of view can not be extracted completely from
the currently existing literature, it will be adopted in this paper,
and we will only deal with class of Poisson structures arising from
Lagrangian morphisms. Thus for us an \emph{\bfseries $n$-shifted
  Poisson structure on $F$} will be defined as an equivalence class of
$n$-shifted Lagrangian maps $F \rightarrow F'$ where $F'$ is a formal
derived stack (in fact in most of the examples we consider $F'$ will
be a derived Artin stack). The typical example is thus the restriction
map
\begin{equation}
  \label{eq:restrictLoc}
Loc_{G}(X) \longrightarrow Loc_{G}(\partial X)
\end{equation}
where $X$ is a compact oriented manifold of dimension $d$ with
boundary $\partial X$. By \cite{ca} this is Lagrangian map and so by
the discussion above it can be considered as a $(2-d)$-shifted Poisson
structure on $Loc_{G}(X)$. When $X$ is a Riemann surface this recovers
the Poisson structure of \cite{fock-rosly,
  ghjw,goldman,guruprasad-rajan}.
In general the bivector underlying the shifted Poisson bracket given
by \eqref{eq:restrictLoc} can be understood explicitly as follows.
The tangent complex of $Loc_{G}(X)$ at a given $G$-local system $\rho$
is $H^*(X,ad(\rho))[1]$. By Lefschetz duality we have a natural
quasi-isomorphism $(H^{\bullet}(X,ad(\rho))[1])^\vee \simeq
H^{\bullet}(X,\partial X,ad(\rho))[d-2]$, and thus a natural element
$$
\xymatrix@1@C+2pc@M+0.5pc{k \ar[r]^-{\mathsf{LD}} &
  H^{\bullet}(X,ad(\rho))[1] \otimes H^{\bullet}(X,\partial X;
  ad(\rho))[d-2].
 }
$$
We can compose this with the boundary map $H^{\bullet}(X,\partial X;
ad(\rho)) \rightarrow 
H^{\bullet}(X,ad(\rho))[1]$
to obtain a map
$$
\xymatrix@C+3pc@R+1pc@M+0.5pc{k \ar[dr]_-{\mathsf{LD}} \ar[r]^-{p} &
  (H^{\bullet}(X,ad(\rho))[1] \otimes H^{\bullet}(X,ad(\rho))[1])[d-2] \\
 &  H^{\bullet}(X,ad(\rho))[1] \otimes 
H^{\bullet}(X,\partial X; ad(\rho))[d-2] \ar[u]
}
$$
This morphism $p$ is the underlying bivector of the $(2-d)$-shifted
Poisson structure on $Loc_{G}(X)$.

Classically a Poisson structure on a smooth variety induces a
foliation of the variety by symplectic leaves.  In our setting, for an
$n$-shifted Poisson structure on a derived stack $F$ given by a
Lagrangian map \linebreak $f : F \rightarrow F'$, the symplectic
leaves are the appropriately interpreted fibers of $f$. Here we need
the qualifier `appropriately interpreted' because we must consider
the fibers in the sense of symplectic geometry, that is  - as fiber
products of Lagrangians in $F'$. Note that specifying a Lagrangian
morphism $\Lambda \rightarrow F'$ is the same thing as specifying a
morphism $* \rightarrow F'$ in the category of Lagrangian
correspondences, and thus is a "point" in this sense. We are therefore
led to the following notion.

\begin{df}\label{d01}
Let $F$ be a derived Artin stack with an $n$-shifted Poisson structure
given by an $n$-shifted Lagrangian morphism $f : F \rightarrow F'$. A
\emph{\bfseries generalized symplectic leaf of $F$} is a derived stack
of the form $F\times_{F'} \Lambda$ for any $n$-shifted Lagrangian
morphism $\Lambda \rightarrow F'$
\end{df}

By \cite{ptvv} a generalized symplectic leaf carries a canonical
$n$-shifted symplectic structure. However, the above definition is a
bit awkward as it depends on the choice of $f$. We will not
try to refine this definition and will take it as a model of
several constructions appearing in the sequel of this work.

Again, the typical example is given by a compact Riemann surface with
boundary $X$.  The restriction map $Loc_{G}(X) \longrightarrow
Loc_{G}(\partial X)=\prod \ [G/G]$ carries a $0$-shifted Lagrangian
structure and thus corresponds to a $0$-shifted Poisson structure on
$Loc_{G}(X)$.  As we have already seen, among the generalized
symplectic leaves of $Loc_{G}(X)$ we have $Loc_{G}(X,\{\lambda_i\})$,
the derived moduli stack of $G$-local systems on $X$ whose monodromies
at infinity are fixed to be conjugate to the given elements $\lambda_i
\in G$.

As a final note, it is instructive to point out that the above notion
of generalized symplectic leaves is a rather flabby notion. For
instance, when the $n$-shifted Poisson structure on $F$ is
non-degenerate (i.e. comes from an $n$-shifted symplectic structure)
then the generalized symplectic leaves are all $n$-shifted symplectic
derived stacks of the form $F \times R$ for some other $n$-shifted
symplectic derived stack $R$. Another example is given by the Poisson
structure on $Loc_{G}(X)$ induced by the restriction map $Loc_{G}(X)
\longrightarrow Loc_{G}(\partial X)$, for an oriented manifold with
boundary. Assume that $Y$ is another oriented manifold with an
identification $\partial Y \simeq \partial X$, then $Loc_{G}(M)$
becomes a generalized symplectic leaf, when $M=Y\sqcup_{\partial
  X}X$. This provides a lot of generalized symplectic leaves, all
given by the different possible ways to complete $X$ to an oriented
manifold without boundary.

\section{Symplectic leaves in the moduli of $G$-local
  systems on smooth complex varieties} \label{sec:4}

In this section we fix a smooth (separated, quasi-compact and
connected) complex algebraic variety $Z$ of complex dimension $d$. We
denote by $X:=Z(\mathbb{C})$ the underlying topological space of
$\mathbb{C}$-points of $X$ endowed with the Eucledian topology. We
also keep the notation $k$ for a given field of characteristic zero
and we fix a reductive group $G$ over $k$ with a chosen non-degenerate
element in $Sym^2(\mathfrak{g}^{\vee})^G$. The derived stack
$Loc_{G}(X)$ is then a derived Artin stack of finite type over $k$ and
we are interested in the following problem:

\begin{pb}\label{pb}
  Show that $Loc_{G}(X)$ carries a natural $(2-2d)$-shifted Poisson
  structure and describe its generalized symplectic leaves.
\end{pb}

\

As we noted before, there are way too many generalized symplectic
leaves according to our definition \ref{d01}. To make the problem more
manageable we will focus on  a  class of generalized symplectic
leaves that is geometrically meaningful. We also want to keep in mind
the case of curves, and when $Z$ is of dimension $1$ we want
our description to recover the symplectic derived stacks
$Loc_{G}(X,\{\lambda_i))$, of $G$-local systems with prescribed
monodromy at infinity.

In the discussion below we will propose a first answer to Problem
\ref{pb}. However, we will restrict ourselves to varieties $Z$ with
nice behavior at infinity. As we will see the problem has a rather
direct and easy answer when the divisor at infinity for $Z$ can be
chosen to be smooth. We will also provide a solution when this divisor
can be chosen to be simple normal crossings with two components where
already some new phenomena arise. We have not analyzed more
complicated behaviors but we are convinced that one can indeed extend
our result to any variety $Z$.

\subsection{The boundary at infinity of a smooth variety}

We start by a general discussion of the notion of boundary at infinity
of a space, and study the specific case of complex algebraic
varieties. These results are not new and we do not claim any
originality, but we record them here for the lack of an adequate
reference.

\begin{df}\label{d1}
The \emph{\bfseries boundary of a topological space $Y$} is by
definition the pro-homotopy type
  $$
\partial Y:=\underset{K \subset Y}{\prolim}(Y-K) \in \op{Pro}(\T),
  $$
where $\prolim$ is the limit taken in the
  $\s$-category $\op{Pro}(\T$) of pro-homotopy types and over the
  opposite category of compact subsets $K \subset Y$.
\end{df}

\

\noindent
The pro-object $\partial Y$ is in general not constant and can be
extremely complicated. However, when $Y= X = Z(\mathbb{C})$ is the
underlying space of a smooth variety $Z$ then $\partial Y$ is
equivalent to a constant pro-object. In fact, more is true:

\begin{prop}\label{p1}
  For a smooth $n$-dimensional complex algebraic variety $Z$ with
  underlying topological space $X = Z(\mathbb{C})$, the pro-object
  $\partial X$ is equivalent to a constant pro-object in $\T$ which
  has the homotopy type of a compact oriented topological manifold of
  dimension $2n-1$.
\end{prop}
{\bfseries Proof:} Let $Z \subset \Zbar$ be a smooth compactification 
such that $D=\Zbar-Z$ is a divisor with simple normal crossing. 
Fix a Riemannian metric on the $C^{\infty}$ manifold underlying 
 $\Zbar$ and  for any $\epsilon >0$ consider the compact subsets
$$
K_{\epsilon}:=\left\{\left. x\in \Zbar  \, \right| \, 
  d(x,D) \geq \epsilon\right\} \subset X.
$$
The system of compact subset $\{K_{\epsilon}\}_{\epsilon \in
  \mathbb{R}_{> 0}}$ is cofinal in the system of all compact subsets
of $X$. We mean here cofinal in the sense of $\s$-category theory, and
the important consequence is that the two pro-objects $\underset{K
  \subset Y}{\prolim}(Y-K)$ and $\underset{\epsilon
  >0}{\prolim}(Y-K_\epsilon)$ are equivalent in the
$\s$-category $\op{Pro}(\T)$.  Moreover, the sets $D_{\epsilon} =
\Zbar - K_{\epsilon}$ of points of distance $< \epsilon$ from $D$
satisfy:
\begin{itemize}
\item For $\epsilon_{1} < \epsilon_{2}$ small enough, the
  inclusion $D_{\epsilon_{1}} \subset D_{\epsilon_{2}}$ is a homotopy
  equivalence;
\item For small enough $\epsilon$ the tubular neighborhood
  $D_{\epsilon}$ retracts to $D$.
\end{itemize}
This is clear near the smooth points of $D$. But near a singular point
$D$ is given by the local equation $z_{1}z_{2} \cdots z_{k} = 0$ for
some local complex analytic coordinates $z_{1}, \ldots, z_{n}$ on
$\Zbar$. In this case the function
$|z_{1}z_{2}\cdots z_{k}|^{2}$ on $\Zbar$ has a non
vanishing gradient and the gradient flow gives the desired retraction
and homotopy equivalence.

Restricting the retraction and homotopy equivalence to the
corresponding punctured tubular neighborhoods $D_{\epsilon} - D =
X - K_{\epsilon}$ we then get that the opens $X -
K_{\epsilon} \subset X$ satisfy:
\begin{itemize}
\item For $\epsilon_{1} < \epsilon_{2}$ small enough the
inclusion $X-K_{\epsilon_{1}} \subset X-K_{\epsilon_{2}}$
is a homotopy equivalence;
\item For $\epsilon$ small enough, $X-K_{\epsilon}$ retracts
  to $\left\{x\in \Zbar \, | \, d(x,D)=\epsilon\right\}$.
\end{itemize}

\

\smallskip

\noindent
This shows that the pro-object $\partial X$ is
equivalent to the constant pro-object $X-K_{\epsilon}$ for
$\epsilon$ small enough and that this constant pro-object is given by 
$\left\{x\in \Zbar \, | \, d(x,D)=\epsilon\right\}$. But
$\left\{x\in \Zbar \, | \, d(x,D)=\epsilon\right\}$ is a
compact submanifold of $X$ of dimension $2n-1$ as this can be
checked locally.  
Indeed if $D$ is given by the local equation
$z_{1}\cdots z_{k} = 0$, then locally the exponenial map on $\Zbar$ gives a
an identification of  $\left\{x\in \Zbar \, | \,
  d(x,D)=\epsilon\right\}$ with the closed subset in
$\mathbb{C}^{n}$ given by the equation $|z_{1}\cdots z_{k}| = \epsilon$.  
It comes, moreover, equiped with a
canonical orientation coming from the complex structure of $X$. 
\hfill $\Box$

\

\begin{rmk} \label{rmk:rob} In the setup of the proof of the previous
  proposition it is instructive to compare the constant pro-object
  $\partial X$ with the boundary of the real oriented blowup of
  $\Zbar$ along the normal crossings divisor $D$. Recall
  \cite{gillam-rob} that given a strict normal crossings divisor $D
  \subset \Zbar$ in a smooth complex algebraic variety, we can form a
  new topological space - the real oriented blowup
  $\op{Bl}_{D}(\Zbar)$ of $\Zbar$ along $D$. The space
  $\op{Bl}_{D}(\Zbar)$ comes with a natural continuous map $\pi :
  \op{Bl}_{D}(\Zbar) \to \Zbar$ and is uniquely characterized (see
  \cite{gillam-rob}) by the properties:
\begin{itemize}
\item[(a)] $\pi : \op{Bl}_{D}(\Zbar)-\pi^{-1}(D) \to \Zbar - D$ is a
  homeomorphism.
\item[(b)] If $(U,z_{1}, \ldots, z_{n})$ is an analytic chart of
  $\Xbar$, such that $U\cap D$ is given by the equation $z_{1}\cdots
  z_{k} = 0$, and if $S^{1} = \{ t \in \mathbb{C} \ | \ |t| = 1 \}$ is
  the unit circle, then 
\[
\pi^{-1}(U) \cong \left\{ (z,t) \in U\times (S^{1})^{k}  \ \left|
\ z_{1} = t_{1} |z_{1}|, \ldots,  z_{k} = t_{k}|z_{k}| \right.\right\},
\]
and in this identification the projection $\pi$ is given by $\pi(z,t) = z$.
\end{itemize}
From (a) and (b) it is clear that $\pi : \op{Bl}_{D}(\Zbar) \to \Zbar$
is defined in the $C^{\infty}$ category and as a $C^{\infty}$ object
$\op{Bl}_{D}(\Zbar)$ is a manifold with corners. As a topological
space $\op{Bl}_{D}(\Zbar)$ is just a topological manifold with
boundary $\delta$ given by the preimage of $D$. The topological
manifold with boundary $\op{Bl}_{D}(\Zbar)$ is homotopy equivalent to
its interior $\op{Bl}_{D}(\Zbar) - \delta = \Zbar - D$ and the pair
$(\op{Bl}_{D}(\Zbar),\delta)$ is homotopy equivalent to the pair
$(X,\partial X)$. Thus $\delta \cong \pi^{-1}(D)$ provides another
model for the constant pro-object $\partial X$.

Note however that the structure of $\delta$ as a $C^{\infty}$ manifold
with corners or even as a stratified topological manifold depends on
the good compactification $\Zbar$ of $X$. Indeed if we replace $\Zbar$
by the usual complex blow up $\widehat{\Zbar}$ of a point $p$ in
$\Zbar$ which is a smooth point in $D$, then $\widehat{\Zbar}$ is a
new good compactification of $X$ whose boundary divisor $\widehat{D}$
has an extra component. The real oriented blow up
$\op{Bl}_{\widehat{D}}(\widehat{\Zbar})$ will have an extra corner and
so will have a boundary $\widehat{\delta}$ which is the same as a
topological manifold but is different as a $c^{\infty}$ manifold with
corners.

This is the reason of why we only view $\partial X$ as the homotopy
type of a topological manifold and not as the isotopy type of a
stratified manifold or a manifold with corners: we need an notion
which is intrinsically associated to $X$, and does not depend on a
particular good compactification.
\end{rmk}

\

\medskip

\noindent
By construction both $X$ and $\partial X$ have the homotopy type of a
finite CW complex, and thus the derived stacks $Loc_{G}(X)$ and
$Loc_{G}(\partial X)$ discussed in the previous section are derived
Artin stacks of finite presentation. Moreover, the canonical map
$\partial X \longrightarrow X$ induces a restriction morphism of
derived Artin stacks
$$
r : Loc_{G}(X) \longrightarrow Loc_{G}(\partial X).
$$ Since the constant pro-object $\partial X$ can be identified with
the topological submanifold \linebreak $\{ x \in \Zbar \ | \ d(x,D) =
\epsilon \}$ of the complex manifold $Z$, we see that $\partial X$
inherits a canonical orientation of dimension $2n-1$. Thus by
\cite{ptvv} the derived stack $Loc_{G}(\partial X)$ carries a
canonical $(3-2n)$-shifted symplectic structure which depends only on
this canonical orientation and on the chosen non-degenerate
$G$-invariant bilinear form on $\mathfrak{g}$.  In fact more is true:
the morphism $\partial X \longrightarrow X$ has the homotopy type of
the inclusion of the boundary of an oriented $2n$-dimensional
manifold.  By \cite{ca} this implies that the restriction morphism
$$
r : Loc_{G}(X) \longrightarrow Loc_{G}(\partial X)
$$
carries a canonical Lagrangian structure with respect to the
canonical shifted symplectic structure on $Loc_{G}(\partial X)$
we just described. On the level of tangent complexes, this
Lagrangian structure reflects Poincar\'e-Lefschetz duality of the manifold
with boundary $(X,\partial X)$. For a given $G$-local system 
$\rho$ on $X$, the Lagrangian structure provides a
natural quasi-isomorphism of complexes
$$
\mathbb{T}_{Loc_{G}(X),\rho} \simeq \mathbb{L}_{Loc_{G}(X)/
  Loc_{G}(\partial X),\rho}[2-2n]
$$
which on  cohomology spaces induces the Poincar\'e
duality isomorphism  on $(X,\partial X)$ with coefficients in $ad(\rho)$:
$$
H^{i}(X,ad(\rho)) \simeq H^{2n-i}(X,\partial X; ad(\rho))^{\vee}.
$$
As explained in the previous section the Lagrangian morphism
$r : Loc_{G}(X) \longrightarrow Loc_{G}(\partial X)$ defines a
canonical $(2-2n)$-shifted Poisson structure on the derived Artin
stack $Loc_{G}(X)$ and so our problem \ref{pb} reduces to the problem
of describing the \emph{\bfseries generalized symplectic leaves} of
$Loc_{G}(X)$.

As we saw in Example~\ref{ex:surface}, when $Z$ is of complex
dimension one the generalized symplectic leaves are obtained by
quasi-Hamiltonian reduction.  Recall that in the language of derived
algebraic geometry the pertinent reductions were constructed as
Lagrangian intersections. Indeed if $\dim_{\mathbb{C}} Z = 1$, the
boundary $\partial X$ has the homotopy type of a disjoint union of
oriented circles, and so the restriction map $r$ can be identified
with a map
$$
r : Loc_{G}(X) \longrightarrow \prod_{i} \ [G/G],
$$
where the product is taken over the points of $\Zbar - Z$ for some
smooth compactification $\Zbar$ of $Z$.  The Lagrangian structure on
this map is equivalent to the data of a quasi-Hamiltonian system, i.e.
of the data of an equivariant group-valued moment map (see
\cite{ca,sa} for details).  Fix elements $\lambda_{i} \in G$ for each
point $i\in \Zbar - Z$, and consider the centralizers $G_{\lambda_{i}}
\subset G$ of the elements $\lambda_{i}$. We have canonical maps
$BG_{\lambda_{i}} \longrightarrow [G/G]$ which are the residual gerbes
of each points $\lambda_{i}$ in $[G/G]$. For the canonical $1$-shifted
symplectic structure on $[G/G]$ each of the maps $BG_{\lambda_{i}}
\longrightarrow [G/G]$ comes equiped with a canonical Lagrangian
structure (for degree reasons the space of Lagrangian structures on
this map is a contractible space). As a result, we can form the
Lagrangian intersection
$$
Loc_{G}(X,\{\lambda_{i}\}):= Loc_{G}(X) \underset{\prod_{i}
  \ [G/G].}{\bigtimes}\prod BG_{\lambda_{i}}
$$
This is the derived Artin stack of $G$-local systems on $Z$ with
local monodromy around the point $i \in \Zbar - Z$ fixed to be in the
conjugacy class of $\lambda_i$. Being a Lagrangian intersection of
$1$-shifted Lagrangian structures this derived stack carries a
canonical $0$-shifted symplectic structure, which on the smooth locus
recovers the well known symplectic structure on symplectic leaves in
character varieties.

Going back to the general case where $Z$ is not necessarly a curve
anymore we again would like to realize the generalized symplectic
leaves of the shifted Poisson derived stack $Loc_{G}(X)$ by an
appropriate quasi-Hamiltonian reduction construction. For this we
start by fixing a good smooth compactification $\Zbar$ of $Z$, i.e. a
smooth proper complex variety $\Zbar$, containing $Z$ as a Zariski
open subset and such that $D=\Zbar-Z$ is a simple normal crossing
divisor.

The idea is to again construct another Lagrangian map
$Loc_{G}(\partial X,\{\lambda_{i}\}) \longrightarrow Loc_{G}(\partial
X)$, where the $\lambda_i$ are elements in $G$ but now $i$ labels the
irreductible components $D_i$ of $D=\Zbar-Z$. In the presence of
intersections of the components of $D$ the construction of the
Lagrangian $Loc_{G}(\partial X,\{\lambda_{i}\}) \rightarrow
Loc_{G}(\partial X)$ appears to be quite complicated.  However, we
analyse below two special cases: the case of a smooth divisor at
infinity and the case where $D$ has only two irreducible components, or 
more generaly has no more than double points (which is enough for the
case of dimension $2$).
We believe that the general case can be handled using similar ideas
but we have not  pursued this direction.

\subsection{The smooth divisor case} 

First we consider the simplest case where $D$ is a smooth divisor - a
disjoint union of connected components $D_i$. In this case $\partial
X$ has the homotopy type of an oriented circle bundle over $\bigsqcup_i
D_i$, which is classified by a collection $\{ \alpha_{i} \}$ where
$\alpha_{i} \in H^{2}(D_{i},\mathbb{Z})$ is the first Chern class of
the normal bundle of $D_i \subset \Zbar$.  Let us fix as above
elements $\lambda_{i} \in G$ with centralizers $G_{\lambda_{i}} \subset
G$.  The group $S^{1}$ acts on the stack $[G/G]=\Map(S^{1},BG)$ by
loop rotations and this action and the cohomology classes $\alpha_{i}$
can be used to define twisted versions
${}_{\alpha_{i}}\widetilde{[G/G]}$ of $[G/G]$ on each $D_{i}$.

To understand this properly we first need to discuss the notion of a
locally constant family of derived stacks over a space. For this, we
recall that any space $T$ can be viewed as a constant derived stack $T
\in \dSt_k$ over $k$. By definition, a \emph{\bfseries family of
  derived stacks over $T$} is a derived stack $F$ together with a map
$F \longrightarrow T$. If $T$ is connected, all the fibers of $F
\longrightarrow T$ are abstractly equivalent as objects in
$\dSt_{k}$. We say that the family has fiber $F_{0}$ if all its fibers
are (non-canonically) equivalent to $F_{0}$. Since $\dSt_{k}$ is an
$\s$-topos, there is an equivalence between $H$-equivariant derived
stacks and derived stacks over $BH$. Below we apply this 
systematically to the case when  $H=S^1=B\mathbb{Z}$.

Write the classes $\alpha_{i}$ as
continuous maps
$$
\alpha_{i} : D_{i} \longrightarrow BS^{1}.
$$
As the group $S^1$ acts on the stack $[G/G]$ we can form the quotient
$[[G/G]/S^1]$ which is a stack over $BS^1$. Using $\alpha_{i}$ we can
pull back $[[G/G]/S^1]$ to $D_{i}$ by $\alpha_i$ to get a locally
constant family of stacks on $D_{i}$, whose fibers are $[G/G]$.  We
denote this family by ${}_{\alpha_{i}}\widetilde{[G/G]} \to D_{i}$ and
we write
$$ \widetilde{[G/G]} = \bigsqcup_{i} \ {}_{\alpha_{i}}\widetilde{[G/G]}
\longrightarrow \bigsqcup_{i} \ D_{i} = D
$$
for the corresponding locally constant family over all of $D$.

Alternatively we can construct $\widetilde{[G/G]}$ as follows. The
class $\alpha_{i}$ defines a circle bundle $\widetilde{D}_{i}
\rightarrow D_{i}$, and so the collection $\{ \alpha_{i} \}$ defines a
circle bundle $p : \widetilde{D} = \bigsqcup_{i} \widetilde{D}_{i} \to
\bigsqcup_{i} D_{i} = D$ over all of $D$. In terms of this projection we
have $\widetilde{[G/G]}\simeq p_*(BG)$ as derived stacks over $D$.

Next observe that for each $i$, the group $S^{1}$ also acts on the
classifying stack $BG_{\lambda_{i}}$, by means of the central element
$\lambda_{i} \in Z(G_{\lambda_{i}}) =
\pi_{1}(\mathsf{aut}(BG_{\lambda_{i}}),\mathsf{id})$. Moreover, for
each $i$ the canonical $1$-shifted Lagrangian map $BG_{\lambda_{i}}
\longrightarrow [G/G]$ comes equipped with a natural
$S^{1}$-equivariant structure for the $S^{1}$ actions on $BG_{\lambda_{i}}$ and
$[G/G]$. Twisting the source and target of this Lagrangian map by
using $\alpha_{i}$ we get locally constant families of stacks
${}_{\alpha_{i}}\widetilde{BG}_{\lambda_{i}} \to D_{i}$ and
${}_{\alpha_{i}}\widetilde{[G/G]} \to D_{i}$, and a $1$-shifted
Lagrangian morphism
\begin{equation} \label{eq:twisted.lag}
  {}_{\alpha_{i}}\widetilde{BG}_{\lambda_{i}} \longrightarrow
  {}_{\alpha_{i}}\widetilde{[G/G]}
\end{equation}
inside the $\s$-category of locally constant families of derived Artin
stacks over $D_{i}$. Since each $D_{i}$ is a compact topological
manifold endowed with a canonical orientation, the map
\eqref{eq:twisted.lag} induces on derived stack of global sections
induces a $(3-2n)$-shifted Lagrangian morphism of derived Artin stacks
$$
r_{i} : Loc_{G_{\lambda_{i}},\alpha_{i}}(D_{i}) \longrightarrow Loc_{G}(\partial_{i} X).
$$
This result is a consequence of the slight generalization of the main
theorem of \cite{ptvv} for which the mapping stacks are replaced by
global sections of locally constant derived stacks. This slight
generalization is proven the exact same way as the case of constant
coefficients, and we will freely use it in this paper.

Here, $\partial_{i} X$ is the connected component of $\partial X$
lying over $D_{i}$, and by definition
$Loc_{G_{\lambda_{i}},\alpha_{i}}(D_{i})$ is the derived stack of
$\alpha_{i}$-twisted principal $G_{\lambda_{i}}$-bundles on $D_{i}$.
  Combining all the $r_{i}$ we get the desired $(3-2n)$-shifted
  Lagrangian morphism
$$
r=\prod_{i} r_{i} : \prod_{i} Loc_{G_{\lambda_{i}},\alpha_{i}}(D_{i}) \longrightarrow
\prod_{i} Loc_{G}(\partial_{i} X) = Loc_{G}(\partial X).
$$
By the Lagrangian intersection theorem of \cite{ptvv} we thus have that
the fiber product of derived stacks
$$
Loc_{G}(X,\{\lambda_{i}\}) := \left(\prod_{i}
Loc_{G_{\lambda_{i}},\alpha_{i}}(D_{i})\right)
\underset{Loc_{G}(\partial X)}{\bigtimes} Loc_{G}(X)
$$
carries a canonical $(2-2n)$-shifted symplectic structure. By
construction/definition, $Loc_{G}(X,\{\lambda_{i}\})$ is the derived
stack of locally constant $G$-bundles on $X$ whose local monodromy
around $D_{i}$ is fixed to be in the conjugacy class
$\mathbf{C}_{\lambda_{i}}$ of $\lambda_{i}$.
Also by construction the  natural projection
$$
Loc_{G}(X,\{\lambda_{i}\}) \longrightarrow Loc_{G}(X)
$$
exhibits
$Loc_{G}(X,\{\lambda_{i}\})$
as a \emph{\bfseries symplectic leaf} of the $(2-2n)$-shifted
Poisson structure on $Loc_{G}(X)$.  

\

\begin{rmk} \label{rem:empty}
Note that the derived stack $Loc_{G_{\lambda_{i},\alpha_{i}}}(D_{i})$
may be empty. Indeed, the groupoid of $k$-points of this stack is the
groupoid of $G$-local systems on $\partial_{i} X$ whose local
monodromy around $D_{i}$ is conjugate to $\lambda_{i}$.  These
$k$-points can also be described as follows. Let $Z(G_{\lambda_{i}})$
be the center of $G_{\lambda_{i}}$. Any
$G_{\lambda_{i}}/Z(G_{\lambda_{i}})$-local system on $D_{i}$
determines a class in $H^{2}(D_{i},Z(G_{\lambda_{i}}))$, which is the
obstruction to lifting this local system to a $G_{\lambda_{i}}$-local
system. For $Loc_{G_{\lambda_{i}},\alpha_{i}}(D_{i})$ to be non-empty
one needs to have a $G_{\lambda_{i}}/Z(G_{\lambda_{i}})$-local system
on $D_{i}$ whose obstruction class matches with the image of
$\alpha_{i}$ under the map $H^{2}(D_{i},\mathbb{Z}) \rightarrow
H^{2}(D_{i},Z(G_{\lambda_{i}}))$ given by $\lambda_{i} : \mathbb{Z}
\to Z(G_{\lambda_{i}})$.  Given $\alpha_{i}$ and $\lambda_{i}$ the
existence of such a local system is a subtle question, closely related
to existence of Azumaya algebras. For instance when $\lambda_{i}$ is a
regular semi-simple element then $G_{\lambda_{i}}$ is a maximal torus
in $G$ (assume $G$ simple and $k$ algebraicaly closed), and thus we
see that the image of $\alpha_i$ in $H^2(D_i,G_{\lambda_{i}})$ must be
zero.  For instance, if in this situation
$\lambda_i$ is of infinte order, this
forces $\alpha_i$ to be a torsion class in $H^2(D_i,\mathbb{Z})$.

\end{rmk}

\subsection{The case of two components} 

We now assume that $D=D_{1} \cup D_{2}$ is the union of two smooth
irreducible components meeting transversally at a smooth codimension
two subvariety $D_{12}=D_1 \cap D_2$. Since the local fundamental
group of $\Zbar - D$ is abelian we fix two commuting elements
$\lambda_{1}$, $\lambda_{2}$ in $G$. Our goal is to construct a
derived moduli stack $Loc_{G}(X,\{\lambda_1,\lambda_{2}\})$ of
$G$-bundles on $X$ with fixed monodromy $\lambda_{1}$ around $D_{1}$
and fixed monodromy $\lambda_{2}$ on around $D_{2}$ and to realize
this stack as a generalized symplectic leaf of $Loc_{G}(X)$.

To set up the problem we need to introduce some notation and
auxiliary stacks. In this setting the homotopy type $\partial X$ can be
represented (see  Remark~\ref{rmk:rob}) as a homotopy push-out
$$
\partial X \simeq \partial_{1} X \bigsqcup_{\partial_{12}X}\partial_{2} X.
$$ Here $\partial_{i} X$ is an oriented circle bundle over
$D_{i}^{o}=D_{i}-D_{12}$, and $\partial_{12}X$ is an oriented $S^1
\times S^1$-bundle over $D_{12}$.  These circle bundles are the
restrictions of the natural circle bundles in
$\mathcal{O}_{\Zbar}(D_{i})$ or equivalently of the natural circle
bundles in the normal bundles of $D_{i}$ in $\Zbar$. The space
$\partial_{12}X$ has the homotopy type of an oriented compact manifold
of dimension $2n-2$, and each component $\partial_i X$ has the
homotopy type of an oriented compact manifold of dimension $2n-1$ with
boundary canonically identified with $\partial_{12}X$. In the same
manner each boundary $\partial (D^{o}_{i})$ is naturally identified
with an oriented $S^{1}$-fibration over $D_{12}$.

For each $D_{i}^{o}$ we have a $\mathbb{Z}$-gerbe on $D_{i}^{o}$ given
by restriction of $\alpha_{i} \in H^{2}(D_{i},\mathbb{Z})$, which is
the restriction of the first Chern class of the normal bundle of
$D_{i}$ inside $\Zbar$. As before, we can form the $\alpha_{i}$-twisted
Lagrangian maps
$$
{}_{\alpha_{i}}\widetilde{BG}_{\lambda_{i}} \longrightarrow
{}_{\alpha_{i}}\widetilde{[G/G]},
 $$
of locally constant derived stacks on $D^{o}_{i}$. We now use the
mapping theorem for manifolds with boundary of \cite{ca} (see also
\cite{ems}) applied to the manifold with boundary $D^{o}_{i}$ and the
Lagrangian map above. Unfolding the definitions we get a Lagrangian
map of derived Artin stacks
$$
\Gamma(D^{o}_{i}; \ {}_{\alpha_{i}}\widetilde{BG}_{\lambda_{i}} )
\longrightarrow \Gamma(\partial(D^{o}_{i});
\ {}_{\alpha_{i}}\widetilde{BG}_{\lambda_{i}})
\underset{\Gamma(\partial(D^{o}_{i});
  \ {}_{\alpha_{i}}\widetilde{[G/G]})}{\bigtimes} \Gamma(D^{o}_{i};
\ {}_{\alpha_{i}}\widetilde{[G/G]}),
$$
where $\Gamma$ here denotes the derived stack of global
sections\footnote{As explained above
derived stacks of global sections are 
the twisted version of derived mapping stacks and can be 
defined formally as being direct images of derived stacks.}. 
By construction, we have
$$
\Gamma(D^{o}_{i}; \ {}_{\alpha_{i}}\widetilde{[G/G]})\simeq
Loc_{G}(\partial_i X) \quad \text{and} \quad
\Gamma(\partial(D^{o}_{i}); \ {}_{\alpha_{i}}\widetilde{[G/G]})\simeq
Loc_G(\partial_{12}X).
$$
We write 
$$
Loc_{G}(\partial_{i} X,\lambda_{i}):=\Gamma(D^{o}_{i}; \
{}_{\alpha_{i}}\widetilde{BG}_{\lambda_{i}}),
$$
for the \emph{\bfseries derived stack of $G$-bundles on $\partial_{i}
  X$ with monodromy $\lambda_{i}$ around $D_{i}^{o}$}.  Similarly we
write
$$
Loc_{G}(\partial_{12}X,\lambda_i):=\Gamma(\partial(D^{o}_{i});
\ {}_{\alpha_{i}} \widetilde{BG}_{\lambda_{i}}),
$$
for the \emph{\bfseries derived stack of $G$-bundles on $\partial_{12}
  X$ with monodromy $\lambda_{i}$ around $D_{i}$}. We can thus rewrite
the above Lagrangian maps as
$$
\ell_i : Loc_{G}(\partial_{i} X,\lambda_{i}) \longrightarrow
Loc_{G}(\partial_{12}X,\lambda_{i})
\underset{Loc_{G}(\partial_{12}X)}{\bigtimes} Loc_{G}(\partial_{i} X).
$$
For $i=1,2$ these are two Lagrangian maps towards an $(3-n)$-shifted
symplectic target. We can consider the direct product, which is still
a Lagrangian morphism
$$
\ell:=\ell_1 \times \ell_2 : Loc_{G}(\partial_{1} X,\lambda_{1})
\times Loc_{G}(\partial_{2} X,\lambda_{2}) \longrightarrow
\prod_{i=1,2}Loc_{G}(\partial_{12}X,\lambda_{i})
\underset{Loc_{G}(\partial_{12}X)}{\bigtimes} Loc_{G}(\partial_{i} X).
$$
Here we think of $\ell$ as a Lagrangian correspondence between two
Lagrangians in $Loc_{G}(\partial_{12}X) \times
Loc_{G}(\partial_{12}X)$, namely
$$
Loc_{G}(\partial_{12} X,\lambda_{1}) \times Loc_{G}(\partial_{12} X,\lambda_{2})
\longrightarrow Loc_{G}(\partial_{12}X) \times 
Loc_{G}(\partial_{12}X)$$
and 
$$
Loc_{G}(\partial_{1} X) \times Loc_{G}(\partial_{2} X) 
\longrightarrow Loc_{G}(\partial_{12}X) \times Loc_{G}(\partial_{12}X).
$$
Pulling back everything to the diagonal of
$Loc_{G}(\partial_{12}X) \times Loc_{G}(\partial_{12}X)$
we get a Lagrangian morphism
$$
\ell : Loc_{G}(\partial_{1} X,\lambda_{1})
\underset{Loc_{G}(\partial_{12}X)}{\bigtimes} Loc_{G}(\partial_{2}
X,\lambda_{2}) \longrightarrow Loc_{G}(\partial X) \times
Loc_{G}(\partial_{12}X,\{\lambda_{1},\lambda_{2}\}),
$$
where we use the short cut notation
$$
Loc_{G}(\partial_{12}X,\{\lambda_{1},\lambda_{2}\})
:=Loc_{G}(\partial_{12}X,\lambda_{1})\underset{
  Loc_{G}(\partial_{12}X)}{\bigtimes} Loc_{G}(\partial_{12}X,\lambda_{2}).
$$
In contrast with the smooth divisor case this setting has an
important new feature, namely the extra term
$Loc_{G}(\partial_{12}X,\{\lambda_{1},\lambda_{2}\})$, which does not
appear when the smooth components of $D$ do not intersect.
Thus, in order to get a Lagrangian map towards $Loc_{G}(\partial X)$ alone
we need to find an extra Lagrangian maping to $Loc_{G}(\partial_{12}
X,\{\lambda_{1},\lambda_{2}\})$. It is not clear to us that such a
Lagrangian always exists, but there is a natural candidate for it that
we will now describe.

We let $G_{(\lambda_{1},\lambda-{2})} = G_{\lambda_{1}}\cap
G_{\lambda_{2}}$ be the centralizer of the pair
$(\lambda_1,\lambda_2)$.  On $D_{12}$, we have a natural
$\mathbb{Z}^{2}$-gerbe, i.e. the external sum $\alpha_{1}\boxplus
\alpha_{2}$ of the restrictions of the two gerbes over $D_{i}$. It
corresponds to $\partial_{12} X$ as a principal $S^{1}\times
S^{1}$-bundle over $D_{12}$. The group $S^{1}\times S^{1}$ acts on the
stack $BG_{(\lambda_{1},\lambda_{2})}$, by the canonical map
$\mathbb{Z}^{2} \rightarrow
\pi_1(\mathsf{aut}(BG_{(\lambda_{1},\lambda_{2})}),
\mathsf{id})=Z(G_{(\lambda_{1},\lambda_{2})})$
given by the pair $(\lambda_1,\lambda_2)$. This provides a twist
${}_{\alpha_{1}\boxplus \alpha_{2}}\widetilde{BG}_{(\lambda_{1},\lambda_{2})}$ of
$BG_{(\lambda_{1},\lambda_{2})}$ on $D_{12}$.  In the same way, we can
define a twist of the double loop stack $\Map(S^{1}\times
S^{1},BG)=[G*G/G]$, where $G*G$ is the derived subscheme of commuting
elements in $G\times G$. That is we have a twisted form
${}_{\alpha_{1}\boxplus \alpha_{2}}\widetilde{[G*G/G]}$ of $[G*G/G]$
over $D_{12}$. The two elements $\lambda_i$, provide natural
inclusion maps
$$
G_{\lambda_{2}} \times \{\lambda_{2}\} \longrightarrow G*G, \qquad
\{\lambda_{1}\} \times G_{\lambda_{1}}
\longrightarrow G*G.
$$ These induce two inclusion maps on quotient stacks
$[G_{\lambda_{i}}/Z_{\lambda_{i}}] \rightarrow [G*G/G]$, which are
naturally $S^1\times S^1$-equivariant.  We thus get maps of twisted
stacks on $D_{12}$, i.e.
${}_{\alpha_{i}}\widetilde{[G_{\lambda_{i}}/G_{\lambda_{i}}]}
\longrightarrow {}_{\alpha_{1}\boxplus\alpha_{2}} \widetilde{[G*G/G]}.
$ Denote the fiber product of these two maps by $\mathcal{F}_{12}$. By
definition, this is a locally constant family of derived Artin stacks
furnished with a fiberwise $(-1)$-shifted symplectic structures, and
so we get an equivalence of derived Artin stacks equipped with
$(3-n)$-shifted symplectic structures
$$
\Gamma(D_{12},\mathcal{F}_{12}) \simeq Loc_{G}(\partial_{12} X,
\{\lambda_1,\lambda_2\}).
$$
There is a canonical point $(\lambda_1,\lambda_2)$ inside
$[G_{\lambda_{1}}/G_{\lambda_{1}}]\times_{[G*G/G]}[G_{\lambda_{2}}/G_{\lambda_{2}}]$
whose stabilizer is $G_{(\lambda_{1},\lambda_{2})}$. This induces a
morphism $BG_{(\lambda_{1},\lambda_{2})} \longrightarrow
[G_{\lambda_{1}}/G_{\lambda_{1}}]\times_{[G*G/G]}[G_{\lambda_{2}}/G_{\lambda_{2}}]$,
which is is $S^1 \times S^1$-equivariant in a natural way.  We
therefore get a twisted version of this map
${}_{\alpha_{1}\boxplus\alpha_{2}}\widetilde{BG}_{(\lambda_{1},\lambda_{2})}
\longrightarrow \mathcal{F}_{12}.$ This map has a canonical isotropic
structure, and to be more precise the space of isotropic structures on
the above map is a contractible space for degree reasons. By taking
global sections we thus obtain an isotropic map
$$
\ell_{12} :
Loc_{G_{(\lambda_{1},\lambda_{2})},\balpha}(D_{12}) \longrightarrow
Loc_{G}(\partial_{12}X,\{\lambda_1,\lambda_2\})),
$$
where $Loc_{G_{(\lambda_{1},\lambda_{2})},\balpha}(D_{12})$ is defined
to be $\Gamma(D_{12};
\ {}_{\alpha_{1}\boxplus\alpha_{2}}\widetilde{BG}_{(\lambda_{1},\lambda_{2})})$.

The question now reduces to understanding whether  the isotropic map
$\ell_{12}$ is Lagrangian. This is the case when the map of derived stacks
$$
BG_{(\lambda_{1},\lambda_{2})} \longrightarrow
[G_{\lambda_{1}}/G_{\lambda_{1}}]\underset{[G*G/G]}{\bigtimes}
[G_{\lambda_{2}}/G_{\lambda_{2}}]
$$
is Lagrangian. A simple examination of the amplitudes of the tangent
complexes shows that this map is Lagrangian if and only if the tangent
complex of
$[G_{\lambda_{1}}/G_{\lambda_{1}}]\times_{[G*G/G]}[G_{\lambda_{2}}/G_{\lambda_{2}}]$
at the canonical point $(\lambda_1,\lambda_2)$ is cohomologically
concentrated in the two extremal degrees $-1$ and $2$. This leads to
the following notion.

\begin{df}\label{d3}
A pair of elements $(\lambda_1,\lambda_2) \in G \times G$ is called
  \emph{\bfseries strict} if it is a commuting pair and if the
  morphism
$$
BG_{(\lambda_{1},\lambda_{2})} \longrightarrow
  [G_{\lambda_{1}}/G_{\lambda_{1}}]\times_{[G*G/G]}[G_{\lambda_{2}}/G_{\lambda_{2}}]
$$
is Lagrangian (for its canonical isotropic structure).
\end{df}

\

\noindent
Assume that $(\lambda_1,\lambda_2)$ is a strict pair. We now have a
new Lagrangian
$$
\ell_{12} :
Loc_{G_{(\lambda_{1},\lambda_{2})},\balpha}(D_{12}) \longrightarrow
Loc_{G}(\partial_{12}X).
$$
By composing with the Lagrangian $\ell$
constructed above, we get the desired Lagrangian map
$$
Loc_{G}(\partial_1 X,\lambda_1)
\underset{Loc_{G}(\partial_{12}X)}{\bigtimes} Loc_{G}(\partial_2
X,\lambda_2)
\underset{Loc_{G}(\partial_{12}X,\{\lambda_1,\lambda_2\}))}{\bigtimes}
Loc_{G_{(\lambda_{1},\lambda_{2})},\alpha}(D_{12}) \longrightarrow
Loc_{G}(\partial X).
$$
The pull-back of this morphism along the restriction map $Loc_{G}(X)
\longrightarrow Loc_{G}(\partial X)$ is thus a derived Artin stack
with a $(2-2n)$-shifted symplectic structure, and its projection to
$Loc_{G}(X)$ can be thought of as a symplectic leaf of the Poisson
structure on $Loc_{G}(X)$.  We denote this symplectic leaf by
$Loc_{G}(X,\{\lambda_1,\lambda_2\})$. We have therefore proven the
following result.

\begin{thm}\label{p3}
With the notation above. 
\begin{enumerate}
\item The derived Artin stack $Loc_{G}(X)$ carries a canonical
  $(2-2n)$-shifted Poisson structure, which is realized by the
  Lagrangian map $Loc_{G}(X) \longrightarrow Loc_{G}(\partial X)$.
\item Let $\Zbar$ be a smooth compactification of $Z$, and assume
  that $\Zbar - Z = D$ is smooth with connected components $D_i$.  Then,
  for any choice of elements $\lambda_i \in G$, the derived Artin stack
  $Loc_{G}(X,\{\lambda_i\})$, of principal $G$-bundles on $X$ whose
  monodromies around $D_i$ in  $\bOO_{\lambda_{i}}$, carries a natural
  $(2-2n)$-shifted symplectic structure and is a symplectic leaf of
  $Loc_{G}(X)$.
\item Let $\Zbar$ be a smooth compactification of $Z$, and assume that 
$\Zbar-Z=D_1\cup D_2$ is a strcit normal crossings divisor  with 
$D_i$ smooth and connected. Then for any commuting pair of elements 
$(\lambda_1,\lambda_2) \in G\times G$ the natural map
  $$
Loc_{G}(\partial_1 X,\lambda_1)
  \underset{Loc_{G}(\partial_{12}X)}{\bigtimes} Loc_{G}(\partial_2
  X,\lambda_2) \longrightarrow Loc_{G}(\partial X) \times
  Loc_{G}(\partial_{12}X,\{\lambda_1,\lambda_2\})
$$
comes equipped with a natural Lagrangian structure.
\item If moreover the pair $(\lambda_{1},\lambda_{2})$  is strict then 
the derived Artin stack
$$Loc_{G}(X,\{\lambda_1,\lambda_2\})$$
comes equipped with a natural $(2-2n)$-shifted symplectic structure which
is a symplectic leaf 
of
$Loc_{G}(X)$.
\end{enumerate}
\end{thm}

\

\begin{rmk} \label{rmk:truncate}
In order to better understand this proposition it is instructive to
examine the situation on the truncated stacks involved. To start with,
the truncation of $Loc_{G}(\partial X)$ is the underived stack of
$G$-local systems on $\partial X$. This can be described as a quotient
stack $[\op{Hom}_{\gp}(\pi_1(\partial X),G)/G]$ (assuming $\partial X$
is connected). The truncation of $Loc_{G}(\partial_1 X,\lambda_1)
\times_{Loc_{G}(\partial_{12}X)} Loc_{G}(\partial_2 X,\lambda_2)$ is
then the full sub-stack consisting of all $G$-local systems on
$\partial X$ for which the local monodromy around $D_i$ is conjugate
to $\lambda_i \in G$.  The truncation of the stack
$Loc_{G}(\partial_{12}X,\{\lambda_1,\lambda_2\})$ is the full
sub-stack of the stack of $G$-local systems on $\partial_{12}X$ whose
local monodromies around $D_i$ are conjugate to $\lambda_i$. Finally,
$Loc_{G_{(\lambda_{1},\lambda_{2})},\balpha}(D_{12})$ is the full sub-stack of
$Loc_{G}(\partial_{12}X,\{\lambda_1,\lambda_2\})$ whose local
monodromy at points of $D_{12}$ is conjugate to the pair
$(\lambda_1,\lambda_2)\in G \times G$.

As a consequence, the truncation of the derived stack
$Loc_{G}(X,\{\lambda_1,\lambda_2\})$ is naturally equivalent to the
full sub-stack of $Loc_{G}(\partial X)$ consisting of $G$-local
systems whose local monodromies around $D_i$ are conjugate to
$\lambda_i$ but also whose local monodromy at points in $D_1 \cap D_2$
is conjugate to the pair $(\lambda_1,\lambda_2)$. Therefore
statement $(4)$ above can be interpretted as the claim that this stack
admits a natural derived structure for which it carries a natural
$(2-2n)$-shifted symplectic structure.
\end{rmk}

\subsection{Strict pairs}

The following result provides many examples
of strict pairs.

\begin{prop}\label{p4}
Let $(\lambda_1,\lambda_2)$ be a commuting pair of elements in $G$,
and $u:=\mathsf{Id}-\mathsf{ad}(\lambda_1)$ and \linebreak
$v:=\mathsf{Id}-\mathsf{ad}(\lambda_2)$ be the corresponding
endormorphisms of $\mathfrak{g}$ induced by the adjoint
representation. Then the pair $(\lambda_1,\lambda_2)$ is strict if and
only $u$ is strict with respect to the kernel of $v$, i.e. we have
$$
\op{Im}(v_{|\ker(u)})=\op{Im}(v)\cap \ker(u).
$$
\end{prop}
{\bfseries Proof:} \ We use the notation introduced above.  Consider
the derived stack
$[G_{\lambda_{1}}/G_{\lambda_{1}}]\times_{[G*G/G]}[G_{\lambda_{2}}/G_{\lambda_{2}}]$
where $G_{\lambda_{i}} \subset G$ is the centralizer of
$\lambda_i$. The derived stack $[G*G/G]=Loc_{G}(S^1\times S^1)$
carries a canonial $0$-shifted symplectic structure, and each map
$[G_{\lambda_{i}}/G_{\lambda_{i}}] \longrightarrow [G*G/G]$ is
Lagrangian. Therefore
$[G_{\lambda_{1}}/G_{\lambda_{1}}]\times_{[G*G/G]}[G_{\lambda_{2}}/G_{\lambda_{2}}]$
carries a canonical $(-1)$-shifted symplectic structure. For degree
reasons the isotropic map
$$ BG_{(\lambda_{1},\lambda_{2})} \longrightarrow
[G_{\lambda_{1}}/G_{\lambda_{1}}]\times_{[G*G/G]}[G_{\lambda_{2}}/G_{\lambda_{2}}]
$$
is Lagrangian if and only if the tangent complex $\mathbb{T}$ of
$[G_{\lambda_{1}}/G_{\lambda_{1}}]\times_{[G*G/G]}
[G_{\lambda_{2}}/G_{\lambda_{2}}]$ taken at the canonical point
$(\lambda_1,\lambda_2)$ is such that
$$
H^0(\mathbb{T})=H^1(\mathbb{T})=0.
$$
As $\mathbb{T}$ is equipped with a $(-1)$-shifted symplectic form, we have 
that $H^0(\mathbb{T})=0$ if and only if $H^1(\mathbb{T})=0$. Therefore, 
the pair $(\lambda_2,\lambda_1)$ is strict if and only if $H^0(\mathbb{T})=0$.

Let $x:=\mathsf{ad}(\lambda_1)$ and $y:=\mathsf{ad}(\lambda_2)$.
The space $H^0(\mathbb{T})$
sits in an five term exact sequence
$$\xymatrix@1{
H^{0}_{y}(H^0_x(\mathfrak{g})) \oplus H^{0}_{x}(H^0_y(\mathfrak{g}))
\ar[r] & H^{0}_{x,y}(\mathfrak{g}) \ar[r]&
H^0(\mathbb{T}) \ar[r] &
H^{1}_{y}(H^0_x(\mathfrak{g})) \oplus H^{1}_{x}(H^0_y(\mathfrak{g}))
\ar[r] &  H^{1}_{x,y}(\mathfrak{g}) \ar[r] & \dots}$$
Here $x$ and $y$ are considered
as  actions of $\mathbb{Z}$ on $\mathfrak{g}$, and 
$H^{\bullet}_x$ and $H^{\bullet}_y$ denote  group cohomology of $\mathbb{Z}$
with coefficients in $\mathfrak{g}$. In the same way $H^{\bullet}_{x,y}$ denotes 
group cohomology of $\mathbb{Z}^2$ with coefficients in $\mathfrak{g}$. 

We have canonical isomorphisms $H^0_x(H^0_y) \simeq H^0_y(H^0_x)
\simeq H^0_{x,y}$ and the first map above is isomorphic to the
sum map on $H^0_{x,y}(\mathfrak{g})$, and therefore is
surjective. This implies that $H^0(\mathbb{T})=0$ if and only if the
last morphism
$$
\phi : H^{1}_{y}(H^0_x(\mathfrak{g})) \oplus
H^{1}_{x}(H^0_y(\mathfrak{g})) \longrightarrow
H^{1}_{x,y}(\mathfrak{g})
$$
is an injective map. Using the Serre spectral sequence for 
the projection to the first factor $\mathbb{Z}^2 \longrightarrow \mathbb{Z}$
we get a short exact sequence
$$
\xymatrix{
0 \ar[r] & H^1_x(H^0_y(\mathfrak{g})) \ar[r] & H^1_{x,y}(\mathfrak{g}) \ar[r]
& H^0_{x}(H^1_y(\mathfrak{g})) \ar[r] & 0.}
$$
The morphism $\phi$ above is compatible with this short exact sequence and
provides a commutative diagram with exact rows
$$\xymatrix{
0 \ar[r] & H^1_x(H^0_y(\mathfrak{g}))  \ar[r] & H^1_{x,y}(\mathfrak{g}) \ar[r]
& H^0_{x}(H^1_y(\mathfrak{g}))  \ar[r] & 0 \\
0 \ar[r] & H^1_x(H^0_y(\mathfrak{g})) \ar[r] \ar[u] & \ar[u]  
H^{1}_{y}(H^0_x(\mathfrak{g})) \oplus H^{1}_{x}
(H^0_y(\mathfrak{g})) \ar[r] 
& H^1_{y}(H^0_x(\mathfrak{g})) \ar[u] \ar[r] & 0.}$$
The map on the left hand side is an identity, and it thus we see that
$H^0(\mathbb{T})=0$ if and only if the natural morphism
$$H^1_{y}(H^0_x(\mathfrak{g})) \longrightarrow H^0_x(H^1_y(\mathfrak{g}))$$
is injective. Unfolding the definition we find the strictness condition
of the proposition. \hfill $\Box$

\

\noindent
Note that since the strictness condition on a pair of elements in $G$
is symmetric by tdefinition, the condition derived in Proposition
\ref{p4} must be symmetric as well.  In particular the role of $u$ and
$v$ in the statement of Proposition \ref{p4} can be exchanged, and so
both conditions are equivalent to each other and equivalent to
strictness. We can use proposition \ref{p4} to produce the following
interesting examples of strict pairs.

\begin{cor}\label{cp4}
Let $(\lambda_1,\lambda_2)$ be a commuting pair of elements in $G$.
\begin{enumerate}
\item If at least one of the $\lambda_i$ is semi-simple 
then the pair $(\lambda_1,\lambda_2)$ is strict.
\item Assume that $\lambda_1$ and $\lambda_2$ are unipotent
elements in $G$, and let 
$$j_i : SL_2 \hookrightarrow G$$
be group embeddings sending $\binom{1\:1}{0\:1}$
to $\lambda_i$. If the
two copies of $SL_2$ in $G$ commute (i.e. 
the $j_1$ and $j_2$ combine into a group homorphism
$j_1\times j_2 : SL_2 \times SL_2 \longrightarrow G$)
then the pair $(\lambda_1,\lambda_2)$ is
  strict.
\end{enumerate}
\end{cor}
{\bfseries Proof:} \ $(1)$ If $\lambda_1$ is semi-simple it defines a
grading on $\mathfrak{g}$ which is preserved by
$v=\mathsf{Id}-\mathsf{a}d(\lambda_2)$. If
$u=\mathsf{Id}-\mathsf{ad}(\lambda_1)$ then $\ker(u)$ is the graded
component of degree $0$, and this obviously implies that strictness
hold $\op{Im}(v_{|\ker(u)})=\op{Im}(v)\cap \ker(u)$.

\

\noindent
$(2)$ The morphism $SL_2 \times SL_2 \longrightarrow G$ 
induces an $SL_2 \times SL_2$-action on the Lie algebra
$\mathfrak{g}$. This action defines  a
decomposition $\mathfrak{g}=\oplus_{p,q}\mathfrak{g}_{p,q}$ of
$\mathfrak{g}$, for which the weights $p$ and $q$ are integers.
Moreover with respect to this decomposition, $u$ acts with bidegree
$(1,0)$ and $v$ acts with bidegree $(0,1)$. Finally, the
Lefschetz property is satisfied:
$$
u_{p,q} : \mathfrak{g}_{p,q} \longrightarrow \mathfrak{g}_{p+1,q}
$$
is injective for $p<0$ and surjective for $p\geq 0$, and similarly
$$
v_{p,q} : \mathfrak{g}_{p,q} \longrightarrow \mathfrak{g}_{p,q+1}
$$ is injective for $q<0$ and surjective for $q\geq 0$. Moreover, 
we have that the map
\[
v_{p,q} : \ker(u_{p,q})
\longrightarrow \ker(u_{p,q+1})
\]
is surjective for all $q\geq 0$.

Let $x \in \op{Im}(v)\cap \ker(u)$. We can decompose
$x=\sum_{p,q}x_{p,q}$ according to the bigrading
$\mathfrak{g}=\oplus_{p,q}\mathfrak{g}_{p,q}$, and by the properties
above we have $x_{p,q}=0$ for $p<0$.  As $x$ lies in the image of $u$,
there are $y_{p,q-1}$ such that $v(y_{p,q-1})=x_{p,q}$. Moreover, for
$q\geq 1$ we can chose $y_{p,q-1} \in \ker(u)$. But if $q\leq 0$, we
have $vu(y_{p,q-1})=u(x_{p,q})=0$, and because $v_{p,q-1}$ is
injective we have $u(y_{p,q-1})=0$. This shows that
$y=\sum_{p,q}y_{p,q}$ is such that $u(y)=0$ and $v(y)=x$. Therefore,
$\op{Im}(v)\cap \ker(u)= \op{Im}(v_{|\ker u})$ and strictness
holds. \hfill $\Box$

\

\begin{rmk} \label{rmk:connections.to.others}
  \ $\bullet$ \, The hypothesis postulating the existence of commuting
  $SL_2$'s in part $(2)$ of the previous corollary is a special case
  of the notion of a Jordan-Lefschetz pair defined and studied by
  Looijenga and Lunts in \cite{looijenga.lunts}.

  \

\noindent  
\ $\bullet$ \, From the proof of part $(2)$ of the above statement, we
see that a stronger result holds: the pair $(\lambda_1,\lambda_2)$ is
strict if a bigrading $\mathfrak{g}=\oplus_{p,q}\mathfrak{g}_{p,q}$ as
in proof exists. Such gradings exist for instance in the setting of
principal nilpotent pairs of \cite{gi}.

\

\noindent
\ $\bullet$ \, Even though the strictness condition is local in nature
there can be global Hodge theoretic reasons that enforce
stricness. For instance the strictness condition will hold
automatically for the local monodromies of any polarized variation of
pure complex Hodge structures on $X$. In fact by the Lefchetz package
for the tame non-abelian Hodge correspondence, the strictness
condition will hold automatically for the local monodromie of any tame
locl system on $X$ which has semismiple global monodromy. We thank
Takuro Mochizuki for this remark.
\end{rmk}

 \

\begin{rmk} \label{rmk:strict.nontrivial}
Finally we note that strictness is a non-trivial condition.  For
instance, if $\lambda$ is any non-trivial unipotent element in $G$,
then the pair $(\lambda,\lambda)$ does not satisfy the strictness
condition of proposition \ref{p4} and thus is not a strict
pair. Indeed in this case $u$ is a non-zero nilpotent endomorphism of
$\mathfrak{g}$ and thus $\ker(u) \cap \op{Im}(u)\neq 0$, but
$\op{Im}(u_{|\ker(u)})=0)$.
\end{rmk}

\subsection{The case of at most double intersection}

The discussion above for a divisor at infinity with at most two smooth
components can be easily extended to the case of any components with
the condition that at most two components intersect at a given
point. This is for instance automatic when $Z$ is a surface.

Assume that we have chosen a compactification $\Zbar$ such that
$D=\Zbar-Z$ can be written a union of smooth connected components
$D=\cup D_i$ for $i=1,\dots,p$. We moreover assume that $D_i\cap D_j$
is connected when non-empty, and we will denote it by $D_{ij}$ (we
always assume $i<j$ here). Finally we assume that $D_i\cap D_j \cap
D_k=\varnothing$ for any three distinct labels $i,j,k$.  As usual we
denote by $D_i^o$ the open in $D_i$ consisting of smooth points of $D$
inside $D_i$.  The boundary $\partial X$ is now (see
Remark~\ref{rmk:rob}) a union of $\partial_i X$ ($S^1$-fibrations over
$D_i^o$) glued together along components of their boundaries
$\partial_{ij}X$ ($S^1\times S^1$-fibrations over $D_{ij}$).

For any $i$ we fix an element $\lambda_i \in G$. We assume that
$(\lambda_i,\lambda_j)$ is a strict pair in the sense of definition
\ref{d3} as soon as $D_{ij}\neq \varnothing$.  Let
$G_{\lambda_{i}}\subset G$ be the centralizer of $\lambda_i$ in $G$.
We have a category $\mathcal{C}_D$, whose objects are the $D_i$ and
the $D_{ij}$ as sub-varieties in $\Zbar$, and whose morphisms are the
inclusions. There is an $\s$-functor
$$
F : \C_D^{op} \longrightarrow \dSt_k
$$
sending each $D_i$ to $Loc_{G}(\partial_i X,\lambda_i)$, the
derived stack of $G$-local systems on $\partial_i X$ whose local
monodromy along $D_i$ is conjugate to $\lambda_i$. By definition the
$\s$-functor $F$ sends $D_{ij}$ to $Loc_{G}(\partial_{ij} X)$, where
$\partial_{ij}X$ is the part of $\partial X$ sitting over $D_{ij}$ as
an $S^1 \times S^1$-bundle. The transition morphisms for the
$\s$-functor $F$ are defined by restriction.

Let $\mathcal{F}$ be the limit of $F$ inside derived stacks. It has a
natural projection to the product
$$
\mathcal{F} \longrightarrow \prod_{i<j}
Loc_{G}(\partial_{ij}X,\{\lambda_i,\lambda_j\}),
$$
where $Loc_{G}(\partial_{ij}X,\{\lambda_i,\lambda_j\})$ is defined as
before.  For each $i<j$ we have a canonical morphism
$$
Loc_{G_{(\lambda_{i},\lambda_{j})},\balpha}(D_{ij}) \longrightarrow
Loc_{G}(\partial_{ij}X,\{\lambda_i,\lambda_j\}),
$$
where
$Loc_{G_{(\lambda_{i},\lambda_{j})},\alpha}(D_{ij})$ is the derived
stack of twisted $G_{(\lambda_{i},\lambda_{j})}$-local systems on
$D_{ij}$ as defined before.  The pull-back possesses a natural
morphism towards $Loc_{G}(\partial X)$
$$
\mathcal{F} \times_{ \prod_{i<j}
  Loc_{G}(\partial_{ij}X,\{\lambda_i,\lambda_j\})} \prod_{i<j}
Loc_{G_{(\lambda_{i},\lambda_{j})},\alpha}(D_{ij}) \longrightarrow
Loc_{G}(\partial X).
$$
This proves the following

\begin{prop}\label{p5}
Under the above assumptions there exists a natural 
Lagrangian structure on the morphism
$$ \mathcal{F} \ \underset{ \prod_{i<j}
  Loc_{G}(\partial_{ij}X,\{\lambda_i,\lambda_j\})}{\bigtimes}
\ \prod_{i<j} Loc_{G_{(\lambda_{i},\lambda_{j})},\balpha}(D_{ij})
\longrightarrow Loc_{G}(\partial X).
$$
\end{prop}

\

\noindent
We can define the derived stack
$Loc_{G}(X,\{\lambda_1,\dots,\lambda_p\})$ as the pull-back of the
Lagrangian in this proposition by the restriction map $Loc_{G}(X)
\rightarrow Loc_{G}(\partial X)$. As a corollary
$Loc_{G}(X,\{\lambda_1,\dots,\lambda_p\})$ carries a natural
$(2-2n)$-shifted symplectic structure. As before, the truncation of
$Loc_{G}(X,\{\lambda_1,\dots,\lambda_p\})$ is the full sub-stack of
$Loc_G(X)$ consisting of all $G$-local systems on $X$ whose local
monodromies around $D_i$ is conjugate to $\lambda_i$, and whose local
monodromies at $D_{ij}$ is conjugate to the strict pair
$(\lambda_i,\lambda_j)$.

\section{Towards a Poisson moduli of connections} \label{sec:5}

We would like to finish this manuscript with some ideas of how to
extend the present results when local systems are replace by bundles
with flat connexions. To start with, for a smooth complex algebraic
variety $X$, it is not possible anymore to use the boundary $\partial
X$, as this would only make sense in the holomorphic
category. Moreover, when $X$ is defined over a smaller field $K
\subset \mathbb{C}$ we also want the moduli of flat bundles on $X$ to
be defined over $K$. As a consequence, if we want to generalize
theorem \ref{p4} to the case of flat bundles a first step is to find
an algebraic counterpart of $\partial X$.

As far as we know there is no algebraic version of $\partial X$,
however several authors have been studying in the recent years a
formal analogue denoted by $\fb X$ (see \cite{bete,ef,hpv}). For a
good compactification $\Zbar$ of $X$, with divisor $D=\Zbar-X$, the
\emph{formal boundary at infinity of $X$} is morally defined as
$\widehat{D}-D$, where $\widehat{D}$ is the formal completion of
$\Zbar$ along $D$. This is only a moral definition as $\widehat{D}-D$
does not actually make sense (it is an empty space when considered in
the sense of formal schemes), but several possible incarnations of
this object have been introduced in \cite{bete,ef,hpv}.  For us, we
follow the approach of \cite{ef} and \cite{hpv}, which do not define
$\fb X$ as an object in its own, but define categories and stacks of
sheaves of perfect complexes $Perf(\fb X)$.  Using the same line of
ideas it is possible to define the derived stack of vector bundles on
$\fb X$ endowed with flat connections $\Vect^{\nabla}(\fb X)$. One key
result, proved in \cite{ef}, is that $\Vect^{\nabla}(\fb X)$ depends
on $X$ alone and not of the chosen compactification $\Zbar$ used to
define it. The derived stack $\Vect^{\nabla}(\fb X)$ is our algebraic
analogue of $Loc_G(\partial X)$ studied in this work. It is then
possible to prove statements analogue to the results mentioned in this
work. As an example we state here a result that will appear in
\cite{pt}.

\begin{thm}\label{tlast}
Let $X$ be a smooth algebraic variety over $k$ of dimension $d$ and
$\Vect^{\nabla}(X)$ the derived stack of vector bundles with flat
connections on $X$.
\begin{enumerate}

\item There is a restriction map $r : \Vect^{\nabla}(X) \rightarrow
  \Vect^{\nabla}(\fb X)$. This map is endowed with a canonical
  Lagrangian structure of degree $2-2d$.

\item The fibers of $r$ are representable by derived quasi-algebraic stacks
  locally of finite presentation.
\end{enumerate}
\end{thm}

Some comments about the previous statement. First of all, we do not
impose any regularity assumption on the connections, and
$\Vect^{\nabla}(X)$ is the derived stack of all connections. In
contrast to the case of local systems, the derived stacks
$\Vect^{\nabla}(X)$ and $\Vect^{\nabla}(\fb X)$ are not representable
as they can have infinite dimensional deformation spaces over general
ring valued points. The meaning of statement $(1)$ is thus subtle as
one has to work with notions such as symplectic and Lagrangian
structures on non-representable objects.  Moreover, the object $\fb X$
does not exist on its own, so the usual constructions methods for
symplectic structures of \cite{ca,ptvv} do not apply as these are
based on evaluation maps which do not exist here.  We overcome this
difficulty by using a completely different construction method, based
on rigid tensor categories and explained in the note \cite{to2}. The
consequence of $(1)$ is of course that $\Vect^{\nabla}(X)$ carries a
canonical Poisson structure.  Finally, the representability statement
$(2)$ states that the derived moduli of flat connections whose formal
structures are fixed at infinity is representable.

We also beleive that symplectic leaves of the Poisson structure on
$\Vect^{\nabla}(X)$ can be defined and studied in a similar fashion as
what we have done in the topological setting. We expect $(2)$ above to
ensure that these symplectic leaves are indeed representable by actual
derived algebraic stacks of finite type. Hopefully, the two results
\ref{p4} and \ref{tlast} can then be related by means of the
Riemann-Hilbert correspondence.  Ultimately one also has to study
derived moduli of Higgs bundles in a similar fashion, and relate to
three kind of moduli spaces by means of the non-abelian Hodge
correspondence of T. Mochizuki.

\smallskip

\noindent
Tony Pantev, {\sc University of Pennsylvania}, tpantev@math.upenn.edu

\smallskip

\noindent
Bertrand To\"{e}n, {\sc  Universit\'e Paul Sabati\'er \& CNRS},
Bertrand.Toen@math.univ-toulouse.fr

\end{document}